\documentclass[a4paper,11pt]{article}
\usepackage[utf8]{inputenc}

\usepackage{amsmath}
\usepackage{amssymb}
\usepackage{amsthm}
\usepackage{esint}
\usepackage{bm}
\usepackage{xcolor}
\usepackage{graphicx}
\usepackage[normalem]{ulem}
\usepackage{upgreek} 
\usepackage{mathtools}
\usepackage{mathrsfs}
\usepackage{subfiles}
\usepackage{enumitem}
\usepackage{hyperref}
\usepackage{accents}
\usepackage[normalem]{ulem}

\newcommand{\el}{\text{el}}
\newcommand{\vi}{\text{vi}}
\newcommand{\eq}{\text{eq}}
\newcommand{\ext}{\mathrm{ext}}

\newcommand{\sym}{\text{sym}}
\newcommand{\an}{\text{anti}}

\newcommand{\R}{\mathbb{R}}
\newcommand{\C}{\mathbb{C}}
\newcommand{\D}{\mathbb{D}}

\newcommand{\M}{\mathbb{M}}

\newcommand{\K}{\mathbb{K}}
\newcommand{\LL}{\mathbb{L}}
\newcommand{\calM}{\mathcal{M}}
\newcommand{\calO}{\mathcal{O}}
\newcommand{\eps}{\varepsilon}
\newcommand{\pl}{\partial}
\newcommand{\DIV}{\mathrm{div\,}}

\newcommand{\Cof}{\mathrm{Cof}}
\newcommand{\scrH}{\mathscr{H}}

\newcommand{\wb}{\overline}
\newcommand{\wh}{\widehat}
\newcommand{\wt}{\widetilde}

\newcommand{\dx}{\mathrm{\,d}x}
\newcommand{\dt}{\mathrm{\,d}t}
\newcommand{\ds}{\mathrm{\,d}s}
\newcommand{\dGamma}{\mathrm{\,d}S}%
\newcommand{\dd}{\mathrm{\,d}}
\newcommand{\dS}{\mathrm{\,d}S}
\newcommand{\id}{\mathrm{id}}

\newcommand{\dist}{\mathrm{dist}}

\newcommand\DT[1]{\mathchoice
                 {{\buildrel{\hspace*{.1em}\text{\LARGE.}}\over{#1}}}
                 {{\buildrel{\hspace*{.1em}\text{\Large.}}\over{#1}}}
                 {{\buildrel{\hspace*{.1em}\text{\large.}}\over{#1}}}
                 {{\buildrel{\hspace*{.1em}\text{\large.}}\over{#1}}}}

\renewcommand{\dot}{\DT}

\newcommand{\calE}{\mathcal{E}}

\newcommand{\calH}{\mathcal{H}}
\newcommand{\calR}{\mathcal{R}}

\newcommand{\rmD}{\mathrm{D}}

\newcommand{\mfh}{\mathfrak{h}}

\newcommand{\SO}{\mathrm{SO}}
\newcommand{\GL}{\mathrm{GL}}
\newcommand{\tr}{\mathrm{tr}}
\DeclarePairedDelimiterX{\abs}[1]{\lvert}{\rvert}{#1}
\DeclarePairedDelimiterX{\Abs}[1]{\big\lvert}{\big\rvert}{#1}
\DeclarePairedDelimiterX{\norm}[1]{\lVert}{\rVert}{#1}
\DeclarePairedDelimiterX{\Norm}[1]{\big\lVert}{\big\rVert}{#1}
\DeclarePairedDelimiterX{\ip}[1]{\langle}{\rangle}{#1}
\newcommand{\wto}{\xrightharpoonup{w}}
\newcommand{\wstarto}{\xrightharpoonup{w^*}}
\newcommand{\sto}{\xrightarrow{s}}

\theoremstyle{plain}
\newtheorem{Th}{Theorem}[section]
\newtheorem*{Th*}{Theorem}
\newtheorem{Lemma}[Th]{Lemma}
\newtheorem{Cor}[Th]{Corollary}

\theoremstyle{definition}
\newtheorem{Def}[Th]{Definition}
\newtheorem*{Def*}{Definition}

\newtheorem{Rem}[Th]{Remark}
\newtheorem{?}[Th]{Problem}
\newtheorem{Ex}[Th]{Example}

\newcommand{\tdots}{\,\vdots\,}

\title{Linearization of finite-strain poro-visco-elasticity\\with degenerate mobility}
\author{Willem J. M. van Oosterhout\footnote{Weierstrass Institute for Applied Analysis and Stochastics, Mohrenstraße 39, Berlin, Germany, \url{willem.vanoosterhout@wias-berlin.de}}}
\date{\today}

\numberwithin{equation}{section}
\usepackage[width=16cm]{geometry}

\begin{document}

\maketitle
\begin{abstract}
A quasistatic nonlinear model for finite-strain poro-visco-elasticity is considered in the Lagrangian frame using Kelvin-Voigt rheology. The model consists of a mechanical equation which is coupled to a diffusion equation with a degenerate mobility. Having shown existence of weak solutions in a previous work, the focus is first on showing boundedness of the concentration using Moser iteration. Afterwards, it is assumed that the external loading is small, and it is rigorously shown that solutions of the nonlinear, finite-strain system converge to solutions of the linear, small-strain system.
\end{abstract}


\vspace{1em}

\section{Introduction}

The theory of elasticity coupled to physical processes, such as diffusion, is a highly relevant topic, both from the theoretical and the applied side. For a non-exhaustive list of applications, we refer to e.g.~\cite{HZZS2008TCDLDPG,CheAna2010CTFPLDEM,WLX*2020CPBT,CaSeSa2022ECMS} for results related to polymeric gels, elastomeric materials, biological tissue and solid-state batteries. For the mathematical analysis of both finite- and small-strain models, we refer to e.g.~\cite{MiRoSa2018GERV,MieRou2020TVKVR,BocWeb2021NQSP,BaFrKr2023NaLMiTVE,vOsLie2023FSPVE}
for results related to elasto-visco-plasticity, thermo-visco-elasticity and poro(-visco)-elasticity.

In the present paper, we are interested in poro-visco-elastic material models in the finite-strain setting, which have been recently investigated concerning their analytical properties in \cite{vOsLie2023FSPVE}. For a time horizon $T>0$ and $\Omega\subset \R^d$ a bounded, open set giving the reference configuration, we look for deformations $\chi:[0,T]\times\Omega\to\R^d$ and concentrations $c:[0,T]\times\Omega\to\R^+$ satisfying the quasi-static system of partial differential equations
\begin{subequations}
		\label{Eqn:Intro:FSSystem}
		\begin{alignat}{2} 
		\label{Eqn:Intro:FS1}
			- \DIV\!\big(\sigma_\el(\nabla\chi,c) + \sigma_\vi(\nabla\chi,\nabla\dot\chi,c) - \DIV \mfh(\rmD^2 \chi)\big) &= f(t)\qquad&\text{in }[0,T]\times\Omega,\\
			\label{Eqn:Intro:FS2}
			\dot{c} - \DIV\!\big(\calM(\nabla\chi,c)\nabla \mu\big) &= 0\qquad&\text{in }[0,T]\times\Omega.
		\end{alignat} 
	\end{subequations}
Here, the total stress $\Sigma_{\text{tot}}=\sigma_\el+\sigma_\vi-\DIV\mfh$ corresponds to a Kelvin-Voigt material, and consists of the elastic stress $\sigma_\el(F,c)=\pl_F\Phi(F,c)$, derived from a free energy density $\Phi(\nabla\chi,c)$, the viscous stress $\sigma_\vi(F,\dot{F},c)=\pl_{\dot{F}}\zeta(F,\dot{F},c)$, coming from a dissipation potential $\zeta(\nabla\chi,\nabla\dot{\chi},c)$, and the higher-order hyperstress $\mfh(G)=\pl_G\scrH(G)$, coming from a potential $\scrH(\rmD^2\chi)$. The function $f$ is a body force density, $\calM$ is the mobility tensor, and $\mu(F,c)=\pl_c\Phi(F,c)$ is the chemical potential.

We remark that the higher-order regularization in the form of the hyperstress $\mfh$ turns the material into a second-grade non-simple material, a notion introduced by Toupin \cite{Toup1962EMwCS}, and used in e.g.~\cite{MieRou2020TVKVR}, \cite{RouTom2020DCEB}. Consequently, it is not necessary to put any convexity assumptions on the free energy density $\Phi$. 

It is important to note that the evolution of the deformation is formulated in the reference configuration, while diffusion processes are usually formulated in the actual configuration. Since our model is formulated completely in the reference configuration, the equation in \eqref{Eqn:Intro:FS2} is therefore a diffusion equation pulled back to the reference confiuration. For example, the Lagrangian mobility tensor $\calM$ is the pull-back of the Eulerian mobility tensor $\M$ via
\begin{equation}
\label{Eqn:MobPB}
\calM(F,c) = \frac{(\Cof F^\top)\M(F,c/\det F)\Cof F}{\det F}\quad \text{for }(F,c)\in\GL^+(d)\times\R^+.
\end{equation}
In contrast to other works such as \cite{Roub2017VMSS}, \cite{RouTom2020DCEB}, \cite{MieRou2020TVKVR}, or \cite{Roub2021CHEC}, we do not assume that the mobility is uniformly positive definite, but instead allow for degenerate mobilities $\calM(F,c)\sim c^m$ (some $m>0$). These mobilities are physically relevant since they model a higher species permeability when the material opens up due to a species concentration increase, see \cite{CheAna2010CTFPLDEM} and Example \ref{Ex:BiotFickDarcy}.

In the previous work \cite{vOsLie2023FSPVE}, it was shown that under suitable assumptions, there always exist weak solutions to the system \eqref{Eqn:Intro:FSSystem} in the sense of Definition \ref{Def:weakSolutionFiniteStrain}. The first result of this paper deals with the regularity of these weak solutions. In particular, we show in Section \ref{Sec:Regularity} using Moser iteration that under slightly stronger assumptions the concentration always stays bounded, which improves the result in \cite{vOsLie2023FSPVE}. 

In the second part of this paper, we assume that the external loading is small, e.g., for $\eps>0$ small, we define $f_*$ by letting $f=:\eps f_*$. We then introduce the rescaled displacement $u_\eps:=\frac{\chi_\eps-\id}{\eps}$ and concentration variation $\rho_\eps:=\frac{c_\eps-c_\eq}{\eps}$, where $c_\eq>0$ is some equilibrium concentration. After rescaling the system \eqref{Eqn:Intro:FSSystem} by $\frac{1}{\eps}$ and letting $\eps\to 0$, we then formally obtain the linearized system
\begin{subequations}
\label{Eqn:Intro:SSSytem}
\begin{alignat}{2} 
  \label{Eqn:Intro:SS1}
- \DIV\!\big(\C e(u) + \K\rho + \D e(\dot{u})
\big) &= f_*(t)\qquad &&\text{in}\ [0,T]\times\Omega,\\
   \label{Eqn:Intro:SS2}
\dot{\rho} - \DIV\!\big(\M(I,c_\eq)\nabla(\K:e(u)+\LL\rho)\big) &= 0\qquad &&\text{in}\ [0,T]\times\Omega.
\end{alignat}
\end{subequations}
Here, we denote by $e(u):=\frac{1}{2}(\nabla u + \nabla u^\top)$ the symmetric part of $\nabla u$, and have introduced the quantities $\C := \pl^2_{FF}\Phi(I,c_\eq)$, $\K := \pl^2_{Fc}\Phi(I,c_\eq)$, $\LL := \pl^2_{cc}\Phi(I,c_\eq)$, and $\D:=\pl^2_{\dot{F}\dot{F}}\zeta(I,0,c_\eq)$. It is important to note that the degeneracy of the mobility $\calM$ disappears in the limit passage as the effective mobility $\M(I,c_\eq)$ is uniformly positive definite. For a physical derivation of these linearized equations using balance laws and thermodynamic principles, we refer to e.g.~\cite{Anand2015DTLP}. The main result in this part of the paper is to show that this limit passage can be done in a rigorous way.

This result is not the first result dealing with the limit passage from finite-strain elasticity to small-strain elasiticy. We refer e.g. to \cite{DMaNePe2002LEaLoFE} for the limit passage in static elasticity, \cite{MieSte2013LPitELoFP} for elasto-plasticity, \cite{FriKru2018PfNTLVE} for visco-elasticity, and \cite{BaFrKr2023NaLMiTVE,BFKM2024PTNT} for thermo-visco-elasticity. 
The novelty in this work is not in the limit passage in the mechanical equation, which is very similar to the limit passage in \cite{FriKru2018PfNTLVE}, \cite{BaFrKr2023NaLMiTVE} and \cite{BFKM2024PTNT}. The novelty, however, is in the limit passage of the diffusion equation, where we have to deal with the degenerate mobility. 
For nondegenerate mobilities, we can test \eqref{Eqn:Intro:FS2} with $\mu$ to obtain a $L^2([0,T]\times\Omega)$-bound for $\nabla\mu$, which can then be used to extract a converging subsequence to pass to the limit. For a degenerate mobility, however, we only get a bound for the flux $\calM\nabla\mu$. Again, we can extract a converging subsequence, but the difficulty is now to identify this limit. 
Furthermore, it should be noted that the linear equations \eqref{Eqn:Intro:SSSytem} are still fully coupled in the sense that both equations depend on the variables $u$ and $\rho$. This is in contrast to \cite{BaFrKr2023NaLMiTVE} and \cite{BFKM2024PTNT}, where, depending on the range of a paramenter, one of the linear equations might only depend on one variable, and be independent of the other.

The paper is structured as follows. In Section \ref{Sec:Setting}, we introduce the model and state the main results. In Section \ref{Sec:Regularity}, we prove the first main result, namely the better regularity for the concentration in the finite-strain setting. Finally, in Section \ref{Sec:LimitPassage} we show the limit passage from the finite-strain poro-visco-elasticity model in \eqref{Eqn:Intro:FSSystem} to the small-strain poro-visco-elasticity model in \eqref{Eqn:Intro:SSSytem}.

\section{Mathematical setting and main results}
\label{Sec:Setting}

\subsection{Notation}
Our model is described in the Lagrangian setting in the 
reference configuration $\Omega\subset \R^d$. We assume that 
$\Omega$ is an open, bounded domain with Lipschitz boundary, and that 
$\pl\Omega = \Gamma_D\cup\Gamma_N$ (disjoint) such that the Dirichlet part has positive surface measures $\int_{\Gamma_D} 1 \dS > 0$. 
We denote by $L^p(\Omega)$, $H^k(\Omega)$, and $W^{k,p}(\Omega)$ the usual Lebesgue and Sobolev spaces with the standard norms, and by $L\log L(\Omega)$ the space of functions $c\in L^1(\Omega)$ for which $\norm{c\log (c)}_{L^1(\Omega)}$ is finite.

We consider deformations $\chi$ on $\Omega$ that are fixed on the Dirichlet part $\Gamma_D$, namely,
we consider the space
\[
W^{2,p}_\id(\Omega;\R^d) := \{\chi\in W^{2,p}(\Omega;\R^d) \mid \chi|_{\Gamma_D} = \id \}.
\] 
Similarly, the (closed) subspace 
$W^{k,p}_0(\Omega)$ denotes the functions in $W^{k,p}(\Omega)$ with zero trace on $\Gamma_D$, e.g.,
\[
H^1_0(\Omega;\R^d) := \{u\in H^1(\Omega;\R^d) \mid u|_{\Gamma_D}=0 \}.
\]

Finally, we denote by $``a\cdot b"$, $``\!A:B"$, and $``G \tdots H "$ 
the scalar products between vectors $a,b\in\R^d$, matrices $A,B\in \R^{d\times d}$, and third-order tensors $G,H\in\R^{d\times d\times d}$, respectively.

\subsection{Finite-strain poro-visco-elasticity}

To model finite-strain poro-visco-elasticity, we denote by $\chi$ the deformation of the material and by $c$ the concentration of some species.
We then consider a free energy density $\Phi=\Phi
(\nabla\chi,c)$, a higher-order regularization $\mathscr{H}
=\mathscr{H}(\rmD^2\chi)$, a dissipation potential $\zeta 
= \zeta(\nabla\chi,\nabla\dot\chi,c)$, and a (Lagrangian) mobility 
tensor $\calM=\calM(\nabla\chi,c)$. The free energy density $\Phi$ 
gives rise to the first Piola--Kirchhoff stress $\sigma_\el$ 
and the chemical potential $\mu$, the
dissipation potential $\zeta$ to the viscous stress $\sigma_\vi$ via
\begin{equation}
\sigma_\text{el}(F,c) := \pl_F\Phi(F,c),\quad 
\mu(F,c) := \pl_c\Phi(F,c),\quad
\text{and}
\quad\sigma_\text{vi}(F,\dot F,c) := \pl_{\dot F}\zeta(F,\dot F, c),
\end{equation}
and the potential $\mathscr{H}$ 
to the hyperstress $\mfh(G) := \pl_G \mathscr{H}(G)$, where we have used 
the placeholders $F$ for $\nabla\chi$, $\dot F$ for $\nabla\dot 
\chi$, and $G$ for $\rmD^2\chi$. 
Following \cite{vOsLie2023FSPVE}, the model coupling the evolution of the deformation $\chi$ and the concentration $c$ is then given in the reference domain $\Omega$ by:
\begin{subequations}
		\label{Eqn:FiniteStrainSystem}
		\begin{align} 
			\label{Eqn:FiniteStrainStress}
			- \DIV\!\big(\sigma_\el(\nabla\chi,c) + \sigma_\vi(\nabla\chi,\nabla\dot\chi,c) - \DIV \mfh(\rmD^2 \chi)\big) &= f(t),\\
			\label{Eqn:FiniteStrainChemPot}
			\dot{c} - \DIV\!\big(\calM(\nabla\chi,c)\nabla \mu\big) &= 0,
		\end{align} 
	\end{subequations}
	completed with the boundary conditions
	\begin{subequations} 
		\label{Eqn:FiniteStrainBoundary}
		\begin{alignat}{2}
			\label{Eqn:FiniteStrainBoundary1}&\chi = \id &\text{on}\ \Gamma_D,\\
			\label{Eqn:FiniteStrainBoundary2}&\big(\sigma_\el(\nabla\chi,c) + \sigma_\vi(\nabla\chi,\nabla\dot\chi,c)\big)\Vec{n} - \DIV_{\!\text{s}}(\mfh(\rmD^2\chi)\Vec{n}) = g(t) \hspace{25pt}&\text{on}\ \Gamma_N,\\
			\label{Eqn:FiniteStrainBoundary3}&\mfh(\rmD^2\chi):(\Vec{n}\otimes\Vec{n}) = 0 &\text{on}\ \pl\Omega,\\
			\label{Eqn:FiniteStrainBoundary4}&\calM(\nabla\chi,c)\nabla\mu\cdot \Vec{n} +\kappa\mu = \kappa \mu_\ext &\text{on}\ \pl\Omega,
		\end{alignat}
	\end{subequations} 
	where $\vec{n}$ denotes the unit normal vector on $\pl\Omega$, and
	$\kappa(x)\geq 0$ and $\mu_{\ext}(t,x)$ are a given permeability and an external potential, respectively. Here, $\DIV\!_\text{s}$ denotes the surface divergence, defined by $\DIV\!_\text{s}(\cdot) = \tr(\nabla_\text{s}(\cdot))$, i.e., the trace of the surface gradient $\nabla_\text{s} v = (I-\Vec{n}\otimes\Vec{n})\nabla v = \nabla v - \frac{\pl v}{\pl\Vec{n}}\Vec{n}$. Finally, we consider initial conditions
	\begin{equation}
		\label{Eqn:FiniteStrainInitial}
		\chi(0) = \chi_0,\quad c(0) = c_0 \quad\text{in}\ \Omega.
	\end{equation}

To prove that the concentration $c$ is bounded, we slightly strengthen the assumptions used in \cite{vOsLie2023FSPVE}, see also Remark \ref{Rem:CompAss} below. Denote for $R>0$ the set
\[
\mathsf{F}_R := \{F\in \GL^+(d)\mid \abs{F}\leq R, \abs{F^{-1}}\leq R, \text{ and }\det F\geq 1/R\}.
\]

\begin{enumerate}[label=\emph{\bfseries (A\arabic*)}]
\item \label{Assu:Hyperstress} The hyperstress potential is a convex, frame-indifferent $C^1$ 
function $\mathscr{H}:\R^{d\times d\times d}\to \R^+$ such 
that the hyperstress is given by $\mfh(G) = \partial_G \mathscr{H}(G)\in\R^{d\times d\times d}$. Moreover, %
there exist $p\in (d,\infty)\cap [3,\infty)$ and constants $C_{\calH,1}$, $C_{\calH,2}$, $C_{\calH,3}>0$ such that 
\[
C_{\calH,1}\abs{G}^p\leq \mathscr{H}(G)\leq C_{\calH,2}(1+\abs{G}^p), \quad\abs{\pl_G\mathscr{H}(G)}\leq C_{\calH,3}\abs{G}^{p-1}\quad \text{for all}\  G\in \R^{d\times d\times d}.
\]

\item \label{Assu:MobilityTensor} 
The mobility tensor $\calM:\GL^+(d)\times \R^+
\to\R^{d\times d}_{\text{sym}}$ is a continuous map.
There exist an exponent $m>0$, and for all $R>0$ there exist constants $C_{0,\calM,R},C_{1,\calM,R}>0$ such that
\begin{equation}
	\begin{gathered}
\xi\cdot \calM(F,c)\xi\geq C_{0,\calM,R} c^{m}\abs{\xi}^2 
\quad\text{and}\quad 		|\calM(F,c)|\leq C_{1,\calM,R} c^{m}\\
\text{for all}\ \xi\in\R^d, ~F\in \mathsf{F}_R, ~c\in \R^+.
	\end{gathered}
\end{equation}
The admissible range of the exponent $m>0$ depends on the 
growth properties of (the derivatives of) $\Phi$ and is fixed in assumption 
\ref{Assu:FreeEn}.

\item \label{Assu:FreeEn} The free energy $\Phi: \GL^+(d)\times\R^+ \to \R$ is bounded from below, continuous, and $C^2$ on $\GL^+(d)\times (0,\infty)$, i.e., for strictly positive concentrations. It is frame indifferent, and satisfies the following assumptions:
\begin{itemize}
    \item[(i)] For any $c\in \R^+$ there exists constants $C_{\Phi,0}, C_{\Phi,1}>0$ and $q\geq\frac{pd}{p-d}$ such that
    \begin{equation*}
\Phi(F,c) \geq C_{\Phi,0}\abs{F} + \frac{C_{\Phi,0}}{(\det F)^q}-C_{\Phi,1}\quad \text{for all}\ F\in \GL^+(d).
\end{equation*}
	\item[(ii)] There exist an exponent $-1<r<\infty$ such that $r+m\geq 0$, and for all $R>0$ constants $C_i: = C_{\Phi,i,R}>0$ ($1\leq i\leq 2$) and constants $\gamma_i := \gamma_{\Phi,i,R}\geq 0$ ($1\leq i\leq 2$) such that 
	\[
	\frac{C_1}{c} + \gamma_1c^r \leq \pl^2_{cc}\Phi(F,c) \leq \frac{C_3}{c} + \gamma_2c^r \quad \text{for all}\ c\in \R^+,~F\in \mathsf{F}_R.
	\]
	Concerning the constants $\gamma_i$, we distinguish two cases: 
 \begin{description}
        \item{\textbf{Case I}:}  We assume $\gamma_1=\gamma_2 = 0$, and also require $1 \leq m\leq 2-\eta$ for some $\eta>0$. 
        \item{\textbf{Case II}:} We assume $\gamma_2\geq\gamma_1>0$. Additionally, we distinguish for this case:
        \begin{itemize}
        	\item{\textbf{Case IIa}:} We require $0<m\leq 3+r-\eta$ for some $\eta>0$.
        \item{\textbf{Case IIb}:} We require $0<m\leq 2-\eta$ for some $\eta>0$.
\end{itemize}         
        
\end{description}

	\item[(iii)] There exist an exponent $\alpha\geq -1$, and for all $R>0$ a constant $C_{\Phi,5,R}>0$ such that 
	\[
	\abs*{\pl^2_{Fc}\Phi(F,c)}\leq C_{\Phi,5,R}c^\alpha \quad \text{for all}\ c\in \R^+,~ F\in\mathsf{F}_R
	\]
	   In \textbf{Case I} above, 
        $\alpha$ is such that $0\leq m+\alpha \leq \frac{p-s}{ps}$, where $1<s=\frac{md+2}{md+1}< 2$ and $0\leq m+2\alpha$.

    In both \textbf{Case IIa} and \textbf{Case IIb} $\alpha$ is such that $0\leq m+\alpha \leq (2+r)\frac{p-s}{ps}$, where $1<s=\min\{\frac{md+2(r+2)}{md+r+2},\frac{d(m+r+1)+2(r+2)}{d(m+r+1)+r+2}\}<2$. Furthermore, in  \textbf{Case IIa} we require that $0\leq m+2\alpha <m+1+r$, while in \textbf{Case IIb} we require that $0\leq m+2\alpha <m+2+2r$

\end{itemize}

\item \label{Assu:EnergyInitFin} For all $R>0$ there exists a concentration $c_R\in \R^+$ such that $\Phi(F,c_R)<\infty$ and $|\pl_c\Phi(F,c_R)|<\infty$ for all $F\in \mathsf{F}_R$.		
	
\item \label{Assu:ViscousStress} 
The dissipation potential
$\zeta:\R^{d\times d}\times \R^{d\times d}\times \R^+\to \R^+$ 
is such that 
$\zeta(F,\dot F,c) = \wh\zeta(\mathscr{C},\dot{\mathscr{C}} ,c)$, where $\mathscr{C} = F^\top F$ is the right Cauchy--Green tensor, and $\dot{\mathscr{C}} = \dot{F}^\top F + F^\top\dot{F}$. 
Here, $\wh\zeta:\R^{d\times d}_\text{sym} \times \R^{d\times d}_\text{sym}\times \R^+ \to \R^+$
is quadratic in the second variable, namely
\[
\wh\zeta(\mathscr{C},\dot{\mathscr{C}},c) 
= \frac{1}{2}\dot {\mathscr{C}}:\mathbb{\wt D}(\mathscr{C},c)\dot {\mathscr{C}}.
\]
We assume that 
there exist constants $C_{\zeta,1}$, $C_{\zeta,2}>0$ such that
the quadratic form fulfills 
\[
C_{\zeta,1}\abs{\dot{\mathscr{C}}}^2 \leq \wh\zeta(\mathscr{C},\dot{\mathscr{C}},c) \leq C_{\zeta,2}\abs{\dot{\mathscr{C}}}^2 \quad \text{for all}\ c\in \R^+,~F\in \R^{d\times d}.
\]

\item \label{Assu:Loading}  The external forces satisfy $f\in W^{1,\infty}(0,T;L^{2}(\Omega;\R^d))$, $g\in W^{1,\infty}(0,T;L^{2}(\pl\Omega;\R^d))$. We set
\[
\ip{\ell(t),\chi} := \int_\Omega f(t)\cdot \chi\dx + 
\int_{\Gamma_N} g(t)\cdot\chi \dS
\] 
such that $\ell\in W^{1,\infty}(0,T;H^{1}(\Omega;\R^d)^*)$.

\item \label{Assu:Permeability} The permeability $\kappa \in L^\infty(\partial\Omega)$ is nonnegative and strictly positive on a part of the boundary $\pl\Omega$ with positive surface measure, i.e., $\int_{\pl\Omega}\kappa\dS \geq \kappa_* > 0$. We assume that the external chemical potential is such that $\mu_\ext\in L^\infty(0,T;L^\infty(\partial\Omega))$.

\item \label{Assu:RegularityInit} The initial conditions satisfy $\chi_0\in W^{2,p}_\id(\Omega;\R^d)$ with $\det\nabla\chi_0\geq a_0>0$ and $c_0\in L^\infty(\Omega)$ with $c_0\geq 0$ and are such that $\int_\Omega \Phi(\nabla\chi_0,c_0)\dd x< \infty$.
\end{enumerate}

\begin{Rem}
\label{Rem:CompAss}
\begin{itemize}
\item[(i)] Compared to the assumptions in \cite{vOsLie2023FSPVE}, we have strengthened some assumptions. First, we have slightly decreased the upper bound of the admissable range of values for $m$, see Assumption \ref{Assu:FreeEn}(ii). Second, we have increased the regularity of $\mu_\ext$ and $c_0$ in Assumption \ref{Assu:Permeability} and \ref{Assu:RegularityInit}, respectively.
\item[(ii)] For the linearization of the equations, we will additionally assume that the external forces are small, see \ref{Assu:ForcesSmall}.
\end{itemize}
\end{Rem}

As main example satisfying these assumptions, we consider the Biot model with linear mobility.

\begin{Ex}[Biot model and Fick/Darcy's law]
	\label{Ex:BiotFickDarcy}
	The Biot model \cite{Biot1941GTo3dC} (see also \cite[Sect.~4]{RouTom2020DCEB})
	with Boltzmann entropy is given by	
	\begin{equation*}
	\Phi(F,c) = \Phi_\el(F) + \frac{1}{2}M_B(c-c_\eq-\beta(\det F{-}1))^2 
	+ k \big(\log\big(\frac{c}{c_\eq}\big)-c+c_\eq\big),
	\end{equation*}
	for some suitable elastic energy $\Phi_\el$ and constants $M_B,\beta,k, c_\eq>0$. In this case, the assumptions for Case II are satisfied with $\alpha=0$ and $r=0$. Defining the (Eulerian) flux as $\boldsymbol{j}=-\M(F,c)\nabla\mu$, and assuming the mobility is linear in $c$, namely $\M(F,c)=c\M_0$ (i.e., $m=1$), we then obtain
	\[
		\boldsymbol{j}=-k\M_0 \nabla c - c\M_0\nabla p, 
	\]
	where $p=M_B(c-c_\eq - \beta(\det F{-}1))$ is the pressure. The first term corresponds to Fick's law, while the second is related 
	to Darcy's law.
\end{Ex}	

We recall the notion of weak solution as introduced in \cite{vOsLie2023FSPVE}.

\begin{Def}[Weak solution finite-strain equations]
		\label{Def:weakSolutionFiniteStrain}
		Let $1<s<2$ be as in \ref{Assu:FreeEn}(iii). We call a pair $(\chi,c)$ a weak solution of the  
		initial-boundary-value problem \eqref{Eqn:FiniteStrainSystem}--\eqref{Eqn:FiniteStrainInitial}
		if $\chi\in L^\infty(0,T;W^{2,p}_\id(\Omega;\R^d))$, $\dot\chi\in L^2(0,T;H^1(\Omega;\R^d))$ and $c\in L^\infty(0,T;L\log L(\Omega))$, $\dot c \in L^s(0,T;W^{1,s}(\Omega)^*)$ with $\nabla c^{\frac{m}{2}}\in L^2(0,T;L^2(\Omega))$ (Case I). In Case II, we additionally require that $c\in L^\infty(0,T;L^{2+r}(\Omega))$ and $\nabla c^{\frac{m+1+r}{2}}$, $\nabla c^{\frac{m}{2}+1+r}\in L^2(0,T;L^2(\Omega))$. The pair $(\chi,c)$ satisfies the integral equations 
		\begin{subequations}
		\label{Eqn:FS:WeakSoln}
        \begin{equation}
		\label{Eqn:FiniteStrain:WeakSoln:chi}
		\int_0^T\int_\Omega \big(\sigma_\el(\nabla\chi,c) + \sigma_\vi(\nabla\chi,\nabla\dot\chi,c)\big) : \nabla \phi 
		+ \mfh(\rmD^2\chi)\tdots \rmD^2\phi\dx\dt 
		= \int_0^T \ip{\ell(t),\phi}\dt
		\end{equation} 
        for all %
		$\phi\in L^2(0,T;W_0^{2,p}(\Omega;\R^d))$, where $\ip{\cdot,\cdot}$ denotes the duality pairing between $W^{2,p}(\Omega;\R^d)^*$ and $W^{2,p}(\Omega;\R^d)$, %
        and
		\begin{equation}
		\label{Eqn:FiniteStrain:WeakSoln:mu}
		\int_0^T \ip{\dot c,\psi}\dt + \int_0^T\int_\Omega\calM(\nabla\chi,c)\nabla\mu \cdot \nabla \psi \dx\dt + \int_0^T\int_{\pl\Omega} \kappa(\mu - \mu_\ext)\psi\dGamma\dt = 0
		\end{equation} 
		for all $\psi\in L^{s'}(0,T;W^{1,s'}(\Omega))$, where $\ip{\cdot,\cdot}$ denotes the duality pairing between $W^{1,s}(\Omega)^*$ and $W^{1,s'}(\Omega)$.
  
		Furthermore, we require that $\mu\in \pl_c\Phi(\nabla\chi,c)$ almost everywhere in $\Omega$, and that $\mu\in L^2([0,T]\times\pl\Omega)$.
              
        \end{subequations}
        
	\end{Def}

It is important to note that $\nabla\mu$ does not exist in the distributional sense in any Lebesgue space. However, using the relation $\nabla\mu = \pl^2_{Fc}\Phi\rmD^2\chi + \pl^2_{cc}\Phi\nabla c$ and the bounds in \ref{Assu:MobilityTensor} and \ref{Assu:FreeEn}, we see that this relation gives a well-defined concept of weak solution. For example, to see that the second integral of \eqref{Eqn:FiniteStrain:WeakSoln:mu} is finite, note that for $F\in \mathsf{F}_R$ we have
\[
\calM(\nabla\chi,c)\nabla\mu \sim c^m(c^\alpha\rmD^2\chi + c^{-1}\nabla c + \gamma c^r\nabla c) \sim c^{m+\alpha}\rmD^2\chi + c^\frac{m}{2}\nabla c^{\frac{m}{2}} + \gamma c^{m+r}\nabla c.
\]
Using the conditions on the exponents in Assumption \ref{Assu:FreeEn} together with H\"older's and the Gagliardo-Nirenberg-Sobolev inequality, it then follows that $\calM(\nabla\chi,c)\nabla\mu\in L^s([0,T]\times\Omega)$, and thus the weak formulation is well-defined.

In \cite{vOsLie2023FSPVE} it was shown that weak solutions exist cf. Definition \ref{Def:weakSolutionFiniteStrain}, i.e, it was proven:
\begin{Th}[{\cite[Thm.~2.7]{vOsLie2023FSPVE}}]
Suppose that the assumptions \ref{Assu:Hyperstress}--\ref{Assu:RegularityInit} hold. Then, the system in \eqref{Eqn:FiniteStrainSystem}--\eqref{Eqn:FiniteStrainInitial} possesses at least one weak solutions in the sense of Definition \ref{Def:weakSolutionFiniteStrain}.
\end{Th}
 
As a first result, we now show that under the slightly strengthened assumptions, we can use Moser iteration to improve the regularity of the concentration. 

\begin{Th}
\label{Th:ConcBounded}
Suppose that the assumptions \ref{Assu:Hyperstress}--\ref{Assu:RegularityInit} hold. Then, any weak solution $(\chi,c)$ of \eqref{Eqn:FiniteStrainSystem}--\eqref{Eqn:FiniteStrainInitial} satisfies $c\in L^\infty(0,T;L^\infty(\Omega))$.
\end{Th}

The proof of this result is postponed to Section \ref{Sec:Regularity}.

\subsection{Passage to small-strain poro-visco-elasticity}

Next, we assume that the external forces are small, and introduce for sufficiently small $\eps>0$ the rescaled forces $f_*$, $g_*$ and $\ell_*$ by $f=:\eps f_*$, $g=:\eps g_*$ and $\ell=:\eps\ell_*$, respectively. Furthermore, we define the rescaled displacement $u_\eps$ and concentration variation $\rho_\eps$ by setting $u_\eps:= \frac{\chi_\eps-\id}{\eps}$ and $\rho_\eps:= \frac{c_\eps-c_\eq}{\eps}$ for some equilibrium concentration $c_\eq>0$. Similarly, the rescaled chemical potential $\mu_{*,\eps}$ is given by $\mu_{*,\eps}:=\frac{\mu_\eps}{\eps}$. The initial conditions $u_0$ and $\rho_0$ are then given as the limits of $u_{0,\eps}=\frac{\chi_{0,\eps}-\id}{\eps}$ and $\rho_{0,\eps}=\frac{c_{0,\eps}-c_\eq}{\eps}$, respectively.

For simplicity, we restrict the analysis to $\kappa\equiv 0$, i.e., the flux satisfies a homogeneous Neumann condition.
We now rewrite the system \eqref{Eqn:FS:WeakSoln} in terms of $u_\eps$ and $\rho_\eps$ and divide both equation by $\eps$, i.e., we now look at the problem
\begin{subequations}
\label{Eqn:FS:WeakScaled}
\begin{multline}
		\label{Eqn:FS:Def:Scaled}
		\int_0^T\int_\Omega \Big(\frac{1}{\eps}\sigma_\el\big(I+\eps\nabla u_\eps,c_\eq+\eps\rho_\eps\big) + \frac{1}{\eps}\sigma_\vi\big(I+\eps\nabla u_\eps,\eps\nabla\dot u_\eps,c_\eq+\eps\rho_\eps\big)\Big) : \nabla \phi \dx\dt\\
		+ \int_0^T\int_\Omega\frac{1}{\eps}\mfh(\eps\rmD^2 u_\eps)\tdots \rmD^2\phi\dx\dt 
		= \int_0^T \ip{\ell_*(t),\phi}\dt
		\end{multline} 
        for
		$\phi\in L^2(0,T;W_0^{2,p}(\Omega;\R^d))$ and
		\begin{equation}
		\label{Eqn:FS:Con:Scaled}
		\int_0^T \ip{\dot \rho_\eps,\psi}\dt + \int_0^T\int_\Omega\calM\big(I+\eps\nabla u_\eps,c_\eq+\eps\rho_\eps\big)\nabla\mu_{*,\eps} \cdot \nabla \psi \dx\dt = 0
		\end{equation} 
		for $\psi\in L^{s'}(0,T;W^{1,s'}(\Omega))$.
\end{subequations}

To pass to the limit $\eps\to 0$ in this system, we now impose the following additional assumptions:
\begin{enumerate}[label=\emph{\bfseries (L\arabic*)}]
\item \label{Assu:HomNeumann} The boundary permeability $\kappa$ satisfies $\kappa\equiv 0$.
\item \label{Assu:PhiNondegen}
There exists a constant $C>0$ such that $\Phi(F,c)\geq C\dist^2(F,\SO(d))$ for all $F\in \GL^+(d), c\in\R^+$, and $\Phi(F,c)=0$ if and only if $F\in SO(d)$ and $c=c_\eq$. 
\item \label{Assu:Regul}
The free energy $\Phi$ is $C^3$ in a neighbourhood of $\SO(d)\times \{c_\eq\}$, and the dissipation potential $\zeta$ is $C^3$ in a neighbourhood of $\SO(d)\times \{0\}\times \{c_\eq\}$.
\item \label{Assu:FreeEnEq}
The material is stress-free when not deformed, i.e., $\sigma_\el(I,c_\eq) = \pl_F\Phi(I,c_\eq) = 0$, and the chemical potential is normalized in the sense that $\mu(I,c_\eq) = \pl_c\Phi(I,c_\eq) = 0$.
\item \label{Assu:MobReg} The mobility tensor $\calM$ is $C^1$ and satisfies for all $R>0$ for some constant $C_R>0$ the bounds
\[
\abs{\rmD_F\calM(F,c)}\leq C_Rc^{m},\quad \abs{\rmD_c\calM(F,c)}\leq C_Rc^{m-1}\quad \text{ for all }F\in\mathsf{F}_R, c\in\R^+.
\]
\item \label{Assu:ForcesSmall} The rescaled external forces satisfy $f_* \in W^{1,\infty}(0,T;L^2(\Omega;\R^d))$ and \\$g_* \in W^{1,\infty}(0,T;L^2(\pl\Omega;\R^d))$.
\item \label{Assu:LinInitRegul}
The initial conditions satisfy $u_0\in W^{2,p}_0(\Omega;\R^d)$ and $\rho_0\in L^\infty(\Omega)$.
\end{enumerate}

\begin{Rem}
We highlight that \ref{Assu:PhiNondegen} requires that the free energy $\Phi$ is always positive, i.e., $\Phi\geq 0$ for all $F\in \GL^+(d)$ and $c\in \R^+$.
\end{Rem}

\subsubsection*{Heuristical derivation of linearized equations}

We now introduce the quantities $\C := \pl^2_{FF}\Phi(I,c_\eq)$, $\K := \pl^2_{Fc}\Phi(I,c_\eq)$, $\LL := \pl^2_{cc}\Phi(I,c_\eq)$, and $\D:=\pl^2_{\dot{F}\dot{F}}\zeta(I,0,c_\eq) = 4\wt\D$. Formally, we now obtain by a Taylor expansion:
\begin{align*}
\sigma_\el\big(I+\eps\nabla u_\eps,c_\eq+\eps\rho_\eps\big) &= \pl_F\Phi\big(I+\eps\nabla u_\eps,c_\eq+\eps\rho_\eps\big)\\
&= \pl_F\Phi(I,c_\eq) + \eps\pl^2_{FF}\Phi(I,c_\eq)\nabla u_\eps + \eps\pl^2_{Fc}\Phi(I,c_\eq)\rho_\eps + \calO(\eps^2)\\
&= \eps(\C\nabla u_\eps + \K\rho_\eps) + \calO(\eps^2),
\end{align*}
\begin{align*}
\sigma_\vi\big(I+\eps\nabla u_\eps,&\nabla\dot u_\eps,c_\eq+\eps\rho_\eps\big) = \pl_{\dot{F}}\zeta\big(I+\eps\nabla u_\eps,\nabla\dot u_\eps,c_\eq+\eps\rho_\eps\big) \\
&= \pl_{\dot{F}}\zeta(I,0,c_\eq) + \eps\pl^2_{\dot{F}F}\zeta(I,0,c_\eq)\nabla u_\eps + \eps\pl^2_{\dot{F}\dot{F}}\zeta(I,0,c_\eq) \nabla\dot u_\eps \\
&\qquad+ \eps\pl^2_{\dot{F}c}\zeta(I,0,c_\eq)\rho_\eps + \calO(\eps^2)\\
&= \eps\D\nabla\dot u_\eps + \calO(\eps^2).
\end{align*}
For the mobility tensor, recall that $\calM$ is the pullback of the Eulerian mobility tensor $\M$ as given in \eqref{Eqn:MobPB}. Then,
\begin{align*}
\calM\big(I+\eps\nabla u_\eps,c_\eq+\eps\rho_\eps\big) &= \frac{(\Cof(I+\eps\nabla u_\eps)^\top)\M(I+\eps\nabla u_\eps,\frac{c_\eq+\eps\rho_\eps}{\det(I+\eps\nabla u_\eps)})\Cof(I+\eps\nabla u_\eps)}{\det(I+\eps\nabla u_\eps)}\\
&= \M(I,c_\eq) + \calO(\eps),
\end{align*}
\begin{align*}
\mu_{*,\eps} &= \pl_{c}\Phi(I+\eps\nabla u_\eps,c_\eq+\eps\rho_\eps)\\
&= \pl_c\Phi(I,c_\eq) + \eps\pl^2_{Fc}\Phi(I,c_\eq):\nabla u_\eps + \eps \pl^2_{cc}\Phi(I,c_\eq)\rho_\eps + \calO(\eps^2)\\
&= \eps(\K:\nabla u_\eps + \LL\rho_\eps) + \calO(\eps^2).
\end{align*}
Using these expansions, we can now take the limit $\eps\to 0$ in \eqref{Eqn:FS:WeakScaled}, and obtain (formally) the linearized equations
\begin{subequations}
\label{Eqn:SS}
\begin{alignat}{1} 
  \label{Eqn:SS:Def}
- \DIV\!\big(\C e(u) + \K\rho + \D e(\dot{u})
\big) &= f_*(t),\\%
   \label{Eqn:SS:Con}
\dot{\rho} - \DIV\!\big(\M(I,c_\eq)\nabla(\K:e(u) + \LL\rho)\big) &= 0,%
\end{alignat}
\end{subequations}
where we denote by $e(u):=\frac{1}{2}(\nabla u + (\nabla u)^\top)$ the symmetric part of $\nabla u$, and define the linearized chemical potential $\mu_*$ by setting $\mu_* := \K:e(u)+\LL\rho$. The fact that $\C$, $\D$, and $\K$ only act on the symmetric part of $\nabla u$ and $\nabla\dot{u}$ follows from the frame-indifference of $\Phi$ and $\zeta$, see Lemma \ref{Lemma:TensorsActOnSymm}. Furthermore, $\C$ is positive definite, see Lemma \ref{Lemma:PropC}.

The boundary and initial conditions now read as follows:
\begin{subequations}
\label{Eqn:SS:Boundary}
		\begin{gather}
			\label{Eqn:SS:Boundary1}u = 0\quad \text{on}\ \Gamma_D,\qquad\big(\C e(u)+ \K\rho+ \D e(\dot u)\big)\Vec{n} = g_*(t)\quad \text{on}\ \Gamma_N,\\
			\label{Eqn:SS:Boundary2}\M(I,c_\eq)\nabla\mu_*\cdot \Vec{n} = 0\qquad \text{on}\ \pl\Omega,
		\end{gather}
	\end{subequations} 
\begin{equation}
		\label{Eqn:SS:Initial}
		u(0) = u_0,\quad \rho(0) = \rho_0 \quad\text{in}\ \Omega.
	\end{equation}

\subsubsection*{Statement of main result}

We now introduce the notion of weak solution for the linearized system \eqref{Eqn:SS}--\eqref{Eqn:SS:Initial}.

\begin{Def}[Weak solution linearized equations]
\label{Def:WeakSolnLin}
A triple $(u,\rho,\mu_*)$ is called a weak solution of the problem \eqref{Eqn:SS}--\eqref{Eqn:SS:Initial} if $u\in L^\infty(0,T;H^1(\Omega;\R^d))$, $\dot u\in L^2(0,T;H^1(\Omega;\R^d))$, $\rho\in L^\infty(0,T;L^2(\Omega))$, $\dot\rho\in L^2(0,T;H^1(\Omega)^*)$, and $\mu_* \in L^2(0,T;H^1(\Omega))$ with $\mu_* = \K e(u) + \LL\rho$. Furthermore, the following integral equations are satisfied
\begin{subequations}
\label{Eqn:SS:WeakSoln}
\begin{equation}
\label{Eqn:SS:WeakDef}
\int_0^T\int_\Omega (\C e(u) + \K\rho + \D e(\dot u)):\nabla\phi \dx\dt 
= \int_0^T \ip{\ell_*(t),\phi}\dt
\end{equation}
for all $\phi\in L^2(0,T;H^1_0(\Omega;\R^d))$, where $\ip{\cdot,\cdot}$ denotes the duality pairing between $H^1(\Omega;\R^d)^*$ and $H^1(\Omega;\R^d)$, and
\begin{equation}
\label{Eqn:SS:WeakCon}
\int_0^T \ip{\dot\rho,\psi}\dt + \int_0^T\int_\Omega \M(I,c_\eq)\nabla\mu_*\cdot \nabla\psi\dx\dt = 0 
\end{equation}
\end{subequations}
for all $\psi\in L^2(0,T;H^1(\Omega))$, where $\ip{\cdot,\cdot}$ denotes the duality pairing between $H^1(\Omega)^*$ and $H^1(\Omega)$. 
\end{Def}

\begin{Ex}[Linear Biot model]
	\label{Ex:LinBiot}
	By setting $\chi=\id+\eps u$, $c=c_\eq+\eps\rho$, and scaling the energy $\Phi$ from Example \ref{Ex:BiotFickDarcy} by $\frac{1}{\eps^2}$, we obtain by letting $\eps\to 0$ the linear Biot model, with free energy
	\begin{equation*}
	\Phi(u,\rho) = \C e(u):e(u) + \frac{1}{2}M_B(\rho-\beta \tr e(u))^2 
	+ \frac{k}{2c_\eq} \rho^2,
	\end{equation*}
	see also e.g.~\cite[Sect.~7.6.1]{KruRoub2019MMICM} or \cite[Sect.~7.1]{Anand2015DTLP}. Note that the degenerate mobility $\M(F,c)=c\M_0$ now reduces to the uniformly positive definite mobility $\M(I,c_\eq)=c_\eq\M_0$.
\end{Ex}	

We now arrive at the main result of the paper.

\begin{Th}[Passage from nonlinear to linear poro-visco-elasticity]
\label{Th:LimitPassage}
Assume that Assumptions \ref{Assu:HomNeumann}--\ref{Assu:LinInitRegul} hold.
\begin{itemize}
\item[(i)]Let $(\chi_\eps, c_\eps)$ be a sequence of weak solutions of the nonlinear system \eqref{Eqn:FS:WeakScaled}. Then, for $u_\eps := \frac{\chi_\eps-\id}{\eps}$, $\rho_\eps := \frac{c_\eps-c_\eq}{\eps}$ and $\mu_{*,\eps} := \frac{\mu_\eps}{\eps}$ we have that (up to subsequences)
\begin{alignat*}{2}
			&u_\eps \wstarto u\quad &&\text{in}\ L^\infty(0,T;H^1(\Omega;\R^d))\cap H^1(0,T;H^1(\Omega;\R^d)),\\
			&\rho_\eps \wstarto \rho\quad&&\text{in}\ L^\infty(0,T;L^{2}(\Omega))\cap H^1(0,T;H^1(\Omega)^*),\\
			&\calM(\nabla\chi_\eps,c_\eps)\nabla\mu_{*,\eps} \wto \M(I,c_\eq)\nabla\mu_*\quad&&\text{in}\ L^2([0,T]\times\Omega),
		\end{alignat*}
		where $(u,\rho)$ is a weak solution of the linearized system \eqref{Eqn:SS}--\eqref{Eqn:SS:Initial}.
\item[(ii)] The weak solution $(u,\rho)$ obtained in (i) is the unique weak solution of \eqref{Eqn:SS}--\eqref{Eqn:SS:Initial}.
\end{itemize} 
\end{Th}

\begin{Rem}
\label{Rem:LinSolns}
\begin{itemize}
\item[(i)] Contrary to the finite-strain setting, in the small-strain setting we have that the gradient of the (linearized) chemical potential $\nabla\mu_*$ is in $L^2([0,T]\times\Omega)$. In particular, note that now $\nabla\rho$ and $\rmD^2 u$ are not in any Lebesgue space. This uses the fact that $\M(I,c_\eq)$ is uniformly positive definite.
\item[(ii)] Since $\rho$ does not denote the concentration (which is always nonnegative), but rather the difference between the concentration $c$ and the equilibrium concentration $c_\eq$, $\rho$ can be both positive and negative.
\end{itemize}
\end{Rem}

\section{Improved regularity for finite-strain solutions}
\label{Sec:Regularity}

We prove that if $(\chi,c)$ is a weak solution of \eqref{Eqn:FiniteStrainSystem}, then we have $c\in L^\infty(0,T;L^\infty(\Omega))$. The proof is based on Moser iteration, see e.g. \cite{Alik1979LpBS}.

\begin{proof}[Proof of Thm.~\ref{Th:ConcBounded}.]
To show the bound $\norm{c}_{L^\infty(0,T;L^\infty(\Omega))} \leq C$, we iterate over an exponent $q$. For brevity, we write $x \lesssim y$ if there exists a constant $C>0$ such that $x\leq Cy$. This constant is independent of $q$, but might depend on $d$, $m$, etc.

\textit{Step 1.}
Let $q\geq 3-m$, then
\begin{align*}
\norm{c(t)}&_{L^q(\Omega)}^q - \norm{c_0}_{L^q(\Omega)}^q =  
\int_0^t\frac{\dd}{\ds}\int_\Omega c(s)^q \dx\ds = q\int_0^t\int_\Omega c^{q-1}\dot c\dx\ds \\
&= -q(q{-}1)\int_0^t\int_\Omega \calM(\nabla\chi,c)c^{q-2}\nabla c\cdot \nabla\mu \dx\ds - q\int_0^t\int_{\pl\Omega}\kappa(\mu - \mu_\ext)c^{q-1}\dS\ds.
\end{align*}

\textit{Step 2.}
Using that $\nabla\mu = \pl^2_{Fc}\Phi\rmD^2\chi + \pl^2_{cc}\Phi\nabla c$ and the bounds on $\Phi$ and $\calM$ in Assumption \ref{Assu:FreeEn} and \ref{Assu:MobilityTensor}, respectively, we then have
\begin{align}
\norm{c(t)}&_{L^q(\Omega)}^q - \norm{c_0}_{L^q(\Omega)}^q \nonumber\\
&\lesssim -q(q{-}1)\int_0^t\int_\Omega \calM(\nabla\chi,c)c^{q-3}\abs{\nabla c}^2\dx\ds \nonumber\\
&\hspace{20pt}-\gamma q(q{-}1)\int_0^t\int_\Omega \calM(\nabla\chi,c)c^{q-2+r}\abs{\nabla c}^2\dx\ds\nonumber\\
&\hspace{20pt}+q(q{-}1)\int_0^t\int_\Omega \calM(\nabla\chi,c)c^{q-2+\alpha}\abs{\nabla c}\abs{\rmD^2\chi}\dx\ds - q\int_0^t\int_{\pl\Omega}\kappa(\mu - \mu_\ext)c^{q-1}\dS\ds\nonumber\\
&\lesssim -q(q{-}1)\int_0^t\int_\Omega c^{q+m-3}\abs{\nabla c}^2\dx\ds \nonumber\\
&\hspace{20pt}-\gamma q(q{-}1)\int_0^t\int_\Omega c^{q+m-2+r}\abs{\nabla c}^2\dx\ds\nonumber\\
&\hspace{20pt}+q(q{-}1)\int_0^t\int_\Omega c^{q+m-2+\alpha}\abs{\nabla c}\abs{\rmD^2\chi}\dx\ds - q\int_0^t\int_{\pl\Omega}\kappa(\mu - \mu_\ext)c^{q-1}\dS\ds\nonumber\\
&\lesssim -\frac{4q(q{-}1)}{(q{+}m{-}1)^2}\Norm{\nabla c^{\frac{q+m-1}{2}}}^2_{L^2([0,t]\times\Omega)} - \gamma_1\frac{4q(q{-}1)}{(q{+}m{+}r)^2}\Norm{\nabla c^{\frac{q+m+r}{2}}}^2_{L^2([0,t]\times\Omega)} \nonumber\\
&\hspace{20pt} + q(q{-}1)\int_0^t\int_\Omega c^{q+m-2+\alpha}\abs{\nabla c}\abs{\rmD^2\chi}\dx\ds - q\int_0^t\int_{\pl\Omega}\kappa(\mu - \mu_\ext)c^{q-1}\dS\ds\nonumber\\ 
&=: I_1 + I_2 + I_3 + I_4.
\end{align}

The first two integrals $I_1$ and $I_2$ are negative and pose no problems. It remains to bound and/or absorb the last two integrals $I_3$ and $I_4$. We distinguish the cases $\gamma_1=\gamma_2=0$ and $\gamma_2\geq\gamma_1>0$ (see Assumption \ref{Assu:FreeEn}(ii)). In particular, recall from Defintion \ref{Def:weakSolutionFiniteStrain} that for the case $\gamma_2\geq\gamma_1>0$ we have slightly better regularity for the concentration $c$.

We first complete the proof for the case $\gamma_1=\gamma_2=0$. Afterwards, we then show the modifications necessary to deal with the case $\gamma_2\geq\gamma_1>0$.

\textit{\underline{Case I: $\gamma_1=\gamma_2=0$.}} %

\textit{Step 3.}
To estimate the third integral $I_3$, we note that with H\"older's inequality and the fact that $c^{q+m-2+\alpha} = c^{\frac{q+m-3}{2}} c^{\frac{q+m-1+2\alpha}{2}}$,
\[
\int_\Omega c^{q+m-2+\alpha}\abs{\nabla c}\abs{\rmD^2\chi}\dx \lesssim \norm{\rmD^2\chi}_{L^p(\Omega;\R^{d\times d\times d})}\Norm{c^{\frac{q+m-3}{2}}\abs{\nabla c}}_{L^2(\Omega)}\Norm{c^{\frac{q+m-1+2\alpha}{2}}}_{L^\frac{2p}{p-2}(\Omega)}.
\]
Thus, it follows that
\begin{align}
\int_\Omega c^{q+m-2+\alpha}\abs{\nabla c}\abs{\rmD^2\chi}\dx
&\lesssim \norm{\rmD^2\chi}_{L^p(\Omega;\R^{d\times d\times d})}\Norm{c^{\frac{q+m-3}{2}}\abs{\nabla c}}_{L^2(\Omega)}\Norm{c^{\frac{q+m-1+2\alpha}{2}}}_{L^\frac{2p}{p-2}(\Omega)}\nonumber\\
&\lesssim \frac{2}{q{+}m{-}1}\Norm{\nabla c^\frac{q+m-1}{2}}_{L^2(\Omega)}\Norm{c^{\frac{q+m-1+2\alpha}{2}}}_{L^\frac{2p}{p-2}(\Omega)}.\label{Eqn:Int1}
\end{align}
Next, we note that by Assumption \ref{Assu:FreeEn}(iii) $-1\leq\alpha<0$, which together with $q\geq 3-m$ implies that $\frac{1}{2}\leq b_1:=\frac{q+m-1+2\alpha}{q+m-1}<1$. 
Thus, we can use the embedding $L^{\frac{2p}{p-2}}(\Omega)\hookrightarrow L^{b_1\frac{2p}{p-2}}(\Omega)$ 
to estimate
\begin{align*}
\Norm{c^{\frac{q+m-1+2\alpha}{2}}}_{L^\frac{2p}{p-2}(\Omega)} = \Norm{c^{\frac{q+m-1}{2}}}^{b_1}_{L^{b_1\frac{2p}{p-2}}(\Omega)}
\lesssim \Norm{c^{\frac{q+m-1}{2}}}^{b_1}_{L^\frac{2p}{p-2}(\Omega)}.
\end{align*}
In particular, since $p>d$ we have $\frac{2p}{p-2}<\frac{2d}{d-2}$, and
we can use the embedding $L^\frac{2d}{d-2}(\Omega)\hookrightarrow L^\frac{2p}{p-2}(\Omega)$ and the Gagliardo--Nirenberg--Sobolev inequality to further estimate the last term to obtain
\begin{align*}
\Norm{c^{\frac{q+m-1+2\alpha}{2}}}^{b_1}_{L^\frac{2p}{p-2}(\Omega)} &\lesssim \Norm{\nabla c^{\frac{q+m-1}{2}}}^{b_1}_{L^2(\Omega)} + \Norm{c^{\frac{q+m-1}{2}}}^{b_1}_{L^1(\Omega)}.
\end{align*}
Combining this with \eqref{Eqn:Int1}, we thus obtain
\begin{align*}
q(q-1)\int_\Omega &c^{q+m-2+\alpha}\abs{\nabla c}\abs{\rmD^2\chi}\dx \\
&\lesssim \frac{2q(q-1)}{q{+}m{-}1} \Norm{\nabla c^{\frac{q+m-1}{2}}}^{1+b_1}_{L^2(\Omega)} + \frac{2q(q-1)}{q{+}m{-}1}\Norm{\nabla c^{\frac{q+m-1}{2}}}_{L^2(\Omega)}\Norm{c^{\frac{q+m-1}{2}}}^{b_1}_{L^1(\Omega)}.
\end{align*}
Since $b_1<1$, we can use Young's inequality with $\epsilon$ to estimate $\Norm{\nabla c^{\frac{q+m-1}{2}}}^{1+b_1}_{L^2(\Omega)} \lesssim \frac{\epsilon}{q}\Norm{\nabla c^{\frac{q+m-1}{2}}}^{2}_{L^2(\Omega)} + C(\epsilon)q^{\frac{1+b_1}{1-b_1}}$ and $\Norm{\nabla c^{\frac{q+m-1}{2}}}_{L^2(\Omega)}\Norm{c^{\frac{q+m-1}{2}}}^{b_1}_{L^1(\Omega)} \lesssim \frac{\epsilon}{q}\Norm{\nabla c^{\frac{q+m-1}{2}}}^{2}_{L^2(\Omega)} + C(\epsilon)q\Norm{c^{\frac{q+m-1}{2}}}^{2b_1}_{L^1(\Omega)}$, and thus
\begin{align*}
q(q-1)\int_\Omega c^{q+m-2+\alpha}\abs{\nabla c}\abs{\rmD^2\chi}\dx &\lesssim q^{\frac{1+b_1}{1-b_1}} + \epsilon\Norm{\nabla c^{\frac{q+m-1}{2}}}^2_{L^2(\Omega)} + q^2\Norm{c^{\frac{q+m-1}{2}}}^{2b_1}_{L^1(\Omega)}\\
&\lesssim q^{\frac{1+b_1}{1-b_1}}+q^2 + \epsilon\Norm{\nabla c^{\frac{q+m-1}{2}}}^2_{L^2(\Omega)} + q^2\Norm{c}^{q+m-1}_{L^{\frac{q+m-1}{2}}(\Omega)}.
\end{align*}
Thus, we conclude that 
\begin{align*}
I_3 &= q(q{-}1)\int_0^t\int_\Omega c^{q+m-2+\alpha}\abs{\nabla c}\abs{\rmD^2\chi}\dx\ds \\
&\lesssim q^{\frac{1+b_1}{1-b_1}}+q^2 + \epsilon\Norm{\nabla c^{\frac{q+m-1}{2}}}^2_{L^2([0,t]\times\Omega)} + q^2\int_0^t\Norm{c}^{q+m-1}_{L^{\frac{q+m-1}{2}}(\Omega)}\ds.
\end{align*}

\textit{Step 4.}
We estimate the boundary integral $I_4$. To do this, we recall that $p>d$ and thus the embedding $W^{2,p}(\Omega)\hookrightarrow C^{1,1-\frac{d}{p}}(\overline{\Omega})$ implies that $\nabla\chi\in L^\infty(0,T;L^\infty(\pl\Omega))$. By Assumption \ref{Assu:FreeEn}(ii) we then have 
\begin{align*}
\mu &= \pl_c\Phi(\nabla\chi,c)-\pl_c\Phi(I,c) + \pl_c\Phi(I,c) -\pl_c\Phi(I,c_\eq) + \pl_c\Phi(I,c_\eq)\\
&\gtrsim - c^\alpha\abs{\nabla\chi-I} + \log\frac{c}{c_\eq} + \pl_c\Phi(I,c_\eq)\\
&\gtrsim \log\frac{c}{c_\eq} - c^\alpha.
\end{align*}
Using also that $\mu_\ext\in L^\infty(0,T;L^\infty(\pl\Omega))$ (see Assumption \ref{Assu:Permeability}), we have
\begin{align*}
- \int_0^t\int_{\pl\Omega}\kappa(\mu - \mu_\ext)c^{q-1}\dS\ds &\lesssim -\int_0^t\int_{\pl\Omega}c^{q-1}\log c\dS\ds + \int_0^t\int_{\pl\Omega} c^{q-1} \dS\ds + 1.
\end{align*}
Since the function $f(x)=-x^{q-1}\log x + Cx^{q-1}$ is bounded from above for $q>1$ and any $C>0$, it follows that
\begin{align*}
I_4 = - q\int_0^t\int_{\pl\Omega}\kappa(\mu - \mu_\ext)c^{q-1}\dS\ds &\lesssim q.
\end{align*}

\textit{Step 5.}
Combining everything, we thus have for $b_0 = m-1$:
\begin{align*}
\norm{c(t)}&_{L^q(\Omega)}^q \lesssim \norm{c_0}_{L^q(\Omega)}^q + q^2\sup_{s\in[0,t]}\Norm{c(s)}^{q+b_0}_{L^\frac{q+m-1}{2}(\Omega)} + q^{\frac{1+b_1}{1-b_1}} + q^2 + q. 
\end{align*}
Now let $q_n := 2^{n}(2-m)+m-1$ and define
\[
a_n := \norm{c_0}_{L^\infty(\Omega)}^{q_n} + \sup_{s\in[0,t]}\norm{c(s)}_{L^{q_n}(\Omega)}^{q_n} + 1.
\]
Then, since $\frac{q_n+m-1}{2}=q_{n-1}$, it follows that for some constants $C_1$, $C_2>0$ we have $a_n\leq C_1(1+C_2^n) a_{n-1}^{(q_n+b_0)/(q_{n-1})}$. Thus, we have 
\begin{align*}
a_n &\leq C_1(1+C_2^n)(C_1(1+C_2^{n-1}))^{\frac{q_n+b_0}{q_{n-1}}}a_{n-2}^{\frac{q_n+b_0}{q_{n-1}}\frac{q_{n-1}+b_0}{q_{n-2}}}\\
&\qquad\vdots\\
&\leq C_1^{1+\frac{q_n+b_0}{q_{n-1}}+\frac{q_n+b_0}{q_{n-1}}\frac{q_{n-1}+b_0}{q_{n-2}}+\ldots+\frac{q_n+b_0}{q_{n-1}}\cdots\frac{q_1+b_0}{q_0}}\\
&\qquad\times (1+C_2^n)(1+C_2^{n-1})^{\frac{q_n+b_0}{q_{n-1}}}\cdots (1+C_2)^{\frac{q_n+b_0}{q_{n-1}}\cdots\frac{q_2+b_0}{q_1}} a_0^{\frac{q_n+b_0}{q_{n-1}}\cdots\frac{q_1+b_0}{q_0}},
\end{align*}
where $q_0=1$.
Now, we note that if $1\leq m\leq 2-\eta$ for some small $\eta>0$, then $\eta 2^i\leq q_i\leq 2^i+1$ and thus
\[
\prod_{i=j}^{n} \frac{q_i+b_0}{q_{i-1}} = \prod_{i=j}^{n} \frac{q_i}{q_{i-1}}\frac{q_i+b_0}{q_i} \leq \prod_{i=j}^n 2\Big(1+\frac{b_0/\eta}{2^i}\Big) \leq 2^{n-j+1}\prod_{i=1}^\infty \Big(1+\frac{b_0/\eta}{2^i}\Big) \leq C2^{n-j+1} \leq Cq_{n-j+1}.
\]
In particular, we obtain
\[
1+\frac{q_n+b_0}{q_{n-1}}+\frac{q_n+b_0}{q_{n-1}}\frac{q_{n-1}+b_0}{q_{n-2}}+\ldots+\frac{q_n+b_0}{q_{n-1}}\cdots\frac{q_1+b_0}{q_0} \leq C\sum_{i=0}^n 2^i = C(2^{n+1}-1) \leq C2^n \leq Cq_n,  
\]
and
\[
n+(n-1)\frac{q_n+b_0}{q_{n-1}}+\ldots+(1)\frac{q_n+b_0}{q_{n-1}}\cdots\frac{q_2+b_0}{q_1} \leq C\sum_{i=1}^n i2^{n-i} = C(2^{n+1}-2-n)\leq C2^n \leq Cq_n.
\]
Thus, after repeated application of Young's inequality, we arrive at 
\[
a_n \leq C_0\big(1+(C_1(1+C_2)a_0)^{Cq_n}\big),
\]
and thus for almost all $t\in[0,T]$:
\[
\norm{c(t)}_{L^{q_n}(\Omega)} \leq Ca_n^{1/q_n} \leq C\big(1+(C_1(1+C_2)a_0)^{Cq_n}\big)^{1/q_n} \leq C\big(\norm{c_0}_{L^\infty(\Omega)} + \sup_{s\in[0,t]}\norm{c(s)}_{L^1(\Omega)} + 1\big)^C.
\] 
Since $m\leq 2-\eta$ for some small $\eta>0$ (see Assumption \ref{Assu:FreeEn}(ii)), we can then take the limit $n\to\infty$ and obtain
\[
\norm{c(t)}_{L^\infty(\Omega)} \leq C\big(\norm{c_0}_{L^\infty(\Omega)} + \sup_{s\in[0,t]}\norm{c(s)}_{L^1(\Omega)} + 1\big)^C\qquad \text{a.e. } t\in[0,T].
\]
Since we have the bound $\norm{c}_{L^\infty(0,T;L^1(\Omega))}\leq C$, this inequality implies the desired bound $\norm{c}_{L^\infty(0,T;L^\infty(\Omega))}\leq C$.
\medskip

\textit{\underline{Case II: $\gamma_2\geq\gamma_1>0$.}} 

\textit{Step 6.}
This case is similar to Case I, but now we can start iterating from the bound $\norm{c}_{L^\infty(0,T;L^{2+r}(\Omega))}\leq C$. 

Note that there are now two possibilities to absorb the integral $I_3$: it is still possible to absorb it in $I_1$, but also in $I_2$. We start by estimating the integral $I_3$ as
\begin{align*}
\int_0^t\int_\Omega c^{q+m-2+\alpha}&\abs{\nabla c}\abs{\rmD^2\chi}\dx\ds 
&\lesssim \frac{2}{q{+}m{+}r}\int_0^t\Norm{\nabla c^\frac{q+m+r}{2}}_{L^2(\Omega)}\Norm{c^{\frac{q+m-2-r+2\alpha}{2}}}_{L^\frac{2p}{p-2}(\Omega)}\ds.
\end{align*}
In Case IIa, we use that $m+2\alpha<m+1+r$ and thus
\[
\Norm{c^{\frac{q+m-2-r+2\alpha}{2}}}_{L^\frac{2p}{p-2}(\Omega)}\lesssim \Norm{c^{\frac{q+m-1}{2}}}^{b_2}_{L^\frac{2p}{p-2}(\Omega)},
\]
where $b_2:= \frac{q+m-2-r+2\alpha}{q+m-1} < 1$.
Using the Gagliardo-Nirenberg-Sobolev inequality on the last term, we continue in the same way as in Case I and find
\begin{align*}
q(q-1)\int_0^t\int_\Omega &c^{q+m-2+\alpha}\abs{\nabla c}\abs{\rmD^2\chi}\dx\ds \\&\lesssim C(q) + \epsilon\Norm{\nabla c^{\frac{q+m-1}{2}}}^2_{L^2([0,t]\times\Omega)} + \epsilon\Norm{\nabla c^{\frac{q+m+r}{2}}}^2_{L^2([0,t]\times\Omega)} + q^2\int_0^t \Norm{c}^{q+m-1}_{L^\frac{q+m-1}{2}(\Omega)}\ds
\end{align*}
for some polynomial function $C(q)$.

To estimate the boundary integral $I_4$, we now note that by Assumption \ref{Assu:FreeEn}(ii) we have 
\begin{align*}
\mu &\gtrsim \log\frac{c}{c_\eq} + c^{r+1}-c_\eq^{r+1} - c^\alpha.
\end{align*}
Thus,
\begin{align*}
- \int_0^t\int_{\pl\Omega}\kappa(\mu - \mu_\ext)c^{q-1}\dS\ds &\lesssim -\int_0^t\int_{\pl\Omega}c^{q-1}\log c\dS\ds - \int_0^t\int_{\pl\Omega} c^{r+q}\dS\ds \\&\qquad+ \int_0^t\int_{\pl\Omega} c^{\alpha+q-1}\dS\ds + \int_0^t\int_{\pl\Omega} c^{q-1} \dS\ds + 1.
\end{align*}
Note that the second integral is nonpositive, and that the third integral can be absorbed in the second using Young's inequality with $\epsilon$. Using that the function $f(x)=-x^{q-1}\log x + Cx^{q-1}$ is bounded from above for $q>1$ and any $C>0$, it then follows that
\begin{align*}
- \int_0^t\int_{\pl\Omega}\kappa(\mu - \mu_\ext)c^{q-1}\dS\ds &\lesssim 1.
\end{align*}
Combining everything, we thus arrive at
\begin{align*}
\norm{c(t)}&_{L^q(\Omega)}^q \lesssim \norm{c_0}_{L^q(\Omega)}^q + q\sup_{s\in[0,t]}\Norm{c(s)}^{q+m-1}_{L^\frac{q+m-1}{2}([0,t]\times\Omega)} + C(q). 
\end{align*}
Setting $q_n := 2^n(3+r-m)+m-1$, we then proceed as before, and finally take the limit $n\to\infty$ and find
\[
\norm{c(t)}_{L^\infty(\Omega)} \leq C\big(\norm{c_0}_{L^\infty(\Omega)} + \sup_{s\in[0,t]}\norm{c(s)}_{L^{2+r}(\Omega)} + 1\big)^C \leq C\qquad \text{a.e. } t\in[0,T].
\]
In Case IIb, we instead use that $m+2\alpha < m+2+2r$ to estimate the integral $I_3$ as
\[
\Norm{c^{\frac{q+m-2-r+2\alpha}{2}}}_{L^\frac{2p}{p-2}(\Omega)}\lesssim \Norm{c^{\frac{q+m+r}{2}}}^{b_3}_{L^\frac{2p}{p-2}(\Omega)},
\]
where $b_3:=\frac{q+m-2-r+2\alpha}{q+m+r}<1$. Again, we can use the Gagliardo-Nirenberg-Sobolev inequality on the last term and proceed as before, and we find
\begin{align*}
\norm{c(t)}&_{L^q(\Omega)}^q \lesssim \norm{c_0}_{L^q(\Omega)}^q + q\sup_{s\in[0,t]}\Norm{c(s)}^{q+m+r}_{L^\frac{q+m+r}{2}([0,t]\times\Omega)} + C(q). 
\end{align*}
Setting $q_n := 2^n(2-m)+m+r$, we then eventually take the limit $n\to\infty$ and again find
\[
\norm{c(t)}_{L^\infty(\Omega)} \leq C\big(\norm{c_0}_{L^\infty(\Omega)} + \sup_{s\in[0,t]}\norm{c(s)}_{L^{2+r}(\Omega)} + 1\big)^C \leq C\qquad \text{a.e. } t\in[0,T].
\]
In both cases we obtain the desired bound $\norm{c}_{L^\infty(0,T;L^\infty(\Omega))}\leq C$, completing the proof.
\end{proof}

\begin{Rem}[Strictly positive lower bound for $c$]

If we assume that \\$\rmD^2\chi\in L^\infty(0,T;L^\infty(\Omega;\R^{d\times d\times d}))$ and that $m-1+2\alpha\geq 0$ (which is only possible in Case II), then we claim it is possible to find a strictly positive lower bound for $c$, i.e., there exists some $c_*>0$ such that $c(t,x)\geq c_*$ for all $(t,x)\in [0,T]\times\Omega$. Indeed, we can set $w=-\min\{0,\log c+K\}$ for sufficiently large $K>0$. Testing \eqref{Eqn:FS:Con:Scaled} with $-q\frac{w_\eps^{q-1}}{c_\eps}$, and using Moser iteration to take the limit $q\to\infty$, we then obtain $\norm{w}_{L^\infty(0,T;L^\infty(\Omega))}\leq C$, which implies the lower bound for $c$. Although the condition on the exponents $m-1+2\alpha\geq 0$ is not too restricting (e.g. the Biot model satisfies this for linear mobility), the condition $\rmD^2\chi\in L^\infty(0,T;L^\infty(\Omega;\R^{d\times d\times d}))$ is in general too restrictive.
\end{Rem}

\section{Limit passage to linearized poro-visco-elasticity}
\label{Sec:LimitPassage}

We start by proving an energy-dissipation inequality for the rescaled solutions $(\chi_\eps,c_\eps)$, which is the time-continuous version of the time-discrete energy-dissipation inequality shown in \cite[Lemma 3.4]{vOsLie2023FSPVE}.

Define the scaled energy $\calE_\eps$ and the dissipation potential $\calR_\eps$ as
\begin{equation}
\calE_\eps(t,u,\rho) := \frac{1}{\eps^2}\int_\Omega \Phi\big(I+\eps\nabla u,c_\eq+\eps\rho\big) + \scrH(\eps\rmD^2u)\dx - \frac{1}{\eps}\ip{\ell_\eps(t),u},
\end{equation}
\begin{equation}
\calR_\eps(u,\dot u,\rho) := \frac{1}{\eps^2}\int_\Omega \zeta(I+\eps\nabla u,\eps\nabla\dot u,c_\eq+\eps\rho)\dx,
\end{equation} 
and the associated quadratic forms by
\begin{equation}
\label{Eqn:E_0}
\calE_0(t,u,\rho) := \int_\Omega \frac{1}{2}\C e(u):e(u) + \rho\K:e(u) + \frac{1}{2}\LL \rho^2 \dx - \ip{\ell_*(t),u},
\end{equation}
\begin{equation}
\label{Eqn:R_0}
\calR_0(\dot u) := \int_\Omega \frac{1}{2}\D e(\dot u):e(\dot u)\dx. 
\end{equation}

\begin{Lemma}[Energy-dissipation inequality] Let $(\chi_\eps,c_\eps)$ be a weak solution of \eqref{Eqn:FS:WeakScaled} in the sense of Definition \ref{Def:weakSolutionFiniteStrain}. Then,
\begin{multline}
\label{Eqn:EnDissBal}
\calE_\eps(t,u_\eps(t),\rho_\eps(t)) + \int_0^t\int_\Omega \calM(I+\eps\nabla u_\eps,c_\eq+\eps\rho_\eps)\nabla\mu_{*,\eps}\cdot\nabla\mu_{*,\eps} \dx\ds + \int_0^t\calR_\eps(u_\eps,\dot u_\eps,\rho_\eps)\ds \\ \leq \calE_\eps(0,u_\eps(0),\rho_\eps(0)) - \int_0^t \ip{\dot\ell_*(t),u_\eps}\ds.
\end{multline}
\end{Lemma}
\begin{proof}
Formally, we obtain the energy-dissipation inequality by testing \eqref{Eqn:FS:Def:Scaled} with $\dot{\chi}_\eps$ and \eqref{Eqn:FS:Con:Scaled} with $\mu_{\eps}$. However, $\mu_{\eps}$ is not a valid test function, see Definition \ref{Def:weakSolutionFiniteStrain}. Instead, we can regularize the diffusion equation by a term $\eta(-\Delta)^\theta\mu$, which was done in \cite{vOsLie2023FSPVE}. Here, $\eta$ denotes the regularization parameter, and $\theta$ is an exponent chosen big enough so that $H^\theta(\Omega)\hookrightarrow L^\infty(\Omega)$. Using this regularization, we can now derive a time-discrete energy-dissipation inequality, which was done in \cite[Lemma~3.4]{vOsLie2023FSPVE}. We can then use lower semicontinuity arguments to obtain the time-continuous energy-dissipation inequality \eqref{Eqn:EnDissBal}. For example, to pass to the limit $(\eta,\tau)\to 0$ in the time-discrete, regularized term $\int_0^{t_k}\int_\Omega \calM(\nabla\wb \chi_{\eta,\tau},\wb c_{\eta,\tau})\nabla\wb\mu_{\eta,\tau}{\cdot}\nabla\wb\mu_{\eta,\tau} \dx\ds$, where a bar denotes the piecewise constant interpolant, and $\tau$ is the time-step, we recall that $\mu=\pl_c\Phi(\nabla\chi,c)$ and thus $\nabla\mu=\pl^2_{Fc}\Phi(\nabla\chi,c)\rmD^2\chi + \pl^2_{cc}\Phi(\nabla\chi,c)\nabla c$. We can then apply \cite[Thm.~7.5]{FonLeo2007MMCV} to the function $F((\nabla\chi,c),(\rmD^2\chi,\nabla c^{\frac{m}{2}})) = \calM(\nabla\chi,c)\nabla\mu(\nabla\chi,c){\cdot}\nabla\mu(\nabla\chi,c)$ and use the strong convergences of $\nabla\wb\chi_{\eta,\tau}$ and $\wb c_{\eta,\tau}$, and the weak convergences of $\rmD^2\wb\chi_{\eta,\tau}$ and $\nabla\wb c_{\eta,\tau}$, which were shown in \cite[Prop.~4.1]{vOsLie2023FSPVE}, to pass to the limit $(\eta,\tau)\to 0$.


%
%
%
%
%
%
%
%
%
%
%
%
%
%
%
%
%
%
%
%
%
%
%
%
%
%
%
%
%
%
%
%
%
%
%
%
%
%

\end{proof}

Before we state the a priori estimates, we recall the following rigidity lemma, taken from \cite[Lemma 4.2]{FriKru2018PfNTLVE}; see also the proof of \cite[Prop. 3.4]{FrRiMu2002GRDN}.
\begin{Lemma}[Rigidity estimates] 
\label{Lemma:RigidityEst}
Let $u\in W^{2,p}_0(\Omega;\R^d)$ be such that $\calE_\eps(t,u,\rho)\leq C$ for all $\rho\in \R$. Then, there exists a $C>0$ such that for a.e. $t\in[0,T]$:
\begin{enumerate}
\item[(i)] $\eps\norm{\nabla u}_{L^2(\Omega)} \leq C\norm{\dist(I+\eps\nabla u,\SO(d))}_{L^2(\Omega)}$, and
\item[(ii)] $\norm{\nabla u}_{L^\infty(\Omega)} \leq C\eps^{-1+\frac{2}{p}}$ .
\end{enumerate}

\end{Lemma}

\subsection{A priori estimates}
We start by proving a priori estimates for the rescaled displacement $u_\eps$.
\begin{Lemma}[A priori estimates, part 1]
\label{Lemma:APriori1}
Let $(\chi_\eps, c_\eps)$ be a weak solution of \eqref{Eqn:FS:WeakScaled} in the sense of Definition \ref{Def:weakSolutionFiniteStrain}. Then, there exists a constant $C>0$ (independent of $\eps$) such that: 
\begin{itemize}
\item[(i)] $\calE_\eps(t,u_\eps,\rho_\eps)\leq C$ for a.e. $t\in[0,T]$,
\item[(ii)] $\norm{u_\eps}_{L^\infty(0,T;H^1(\Omega;\R^{d}))} \leq C$,
\item[(iii)] $\norm{\nabla\dot u_\eps}_{L^2([0,T]\times\Omega)}\leq C$,
\item[(iv)] $\norm{\rmD^2u_\eps}_{L^\infty(0,T;L^p(\Omega))}\leq C\eps^{-1+\frac{2}{p}}$.
\end{itemize}
\end{Lemma}

\begin{proof}
\begin{itemize}
\item[(i)] From the energy-dissipation inequality \eqref{Eqn:EnDissBal} it follows that for \\$\Lambda
=\norm{\dot\ell_*}_{L^\infty(0,T;W^{2,p}(\Omega;\R^d)^*)}$ we have
\[
\calE_\eps(t,u_\eps,\rho_\eps)\leq C\Big(1+\frac{\Lambda}{\eps^2}\int_0^t\norm{\chi_\eps}_{W^{2,p}(\Omega;\R^d)}\dt\Big).
\]
On the other hand, by coercivity of $\Phi$ (see Assumption \ref{Assu:FreeEn}(i)):
\[
\calE_\eps(t,\chi_\eps,c_\eps) \geq \frac{1}{\eps^2}\Big(C_1\norm{\chi_\eps}_{W^{2,p}(\Omega;\R^d))} + C_2\norm{(\det\nabla\chi_\eps)^{-1}}^q_{L^q(\Omega))} - C_3\Big).
\]
So, using Gr\"onwall's lemma, the bound follows.

\item[(ii)] 
Using the rigidity estimate $\eps\norm{\nabla u_\eps}_{L^2(\Omega;\R^{d\times d})}\leq C\norm{\dist(\nabla\chi_\eps,\SO(d)}_{L^2(\Omega; \R^{d\times d})}$ from Lemma \ref{Lemma:RigidityEst}(i), Assumption \ref{Assu:PhiNondegen}, and the boundedness of $\calE_\eps$ it follows that 
\begin{align*}
\norm{\nabla u_\eps}^2_{L^\infty(0,T;L^2(\Omega;\R^{d\times d}))} &\leq \frac{1}{\eps^2}\sup_{t\in[0,T]}\int_\Omega \dist^2(\nabla\chi_\eps,\SO(d))\dx \\
&\leq \frac{1}{\eps^2}\sup_{t\in[0,T]}\int_\Omega\Phi(\nabla\chi_\eps,c_\eps)\dx \leq C.
\end{align*}
The bound $\norm{u_\eps}_{L^\infty(0,T;H^1(\Omega;\R^{d}))} \leq C$ now follows from Poincar\'e's inequality.

\item[(iii)] Since $\nabla\chi_\eps\in \mathsf{F}_R$ for some $R>0$, we can use the generalized Korn's inequality as in \cite[Cor.\ 3.4]{MieRou2020TVKVR} and Assumption \ref{Assu:ViscousStress} to obtain
		\begin{align*}
			\norm{\dot u_\eps}^2_{L^2(0,T;H^1(\Omega))} &= 
			\frac{1}{\eps^2}\norm{\dot\chi_\eps}^2_{L^2(0,T;H^1(\Omega))}
			= \frac{1}{\eps^2}\int_0^T \norm{\dot\chi_\eps(t)}^2_{H^1(\Omega)}\dt\\
			&\leq \frac{C}{\eps^2}\int_0^T\int_\Omega \abs*{(\nabla\chi_\eps)^\top\nabla\dot\chi_\eps + (\nabla\dot\chi_\eps)^\top\nabla\chi_\eps}^2\dx\dt\\
			&\leq \frac{C}{\eps^2}\int_0^T\int_\Omega \wh\zeta\big((\nabla\chi_\eps)^\top\nabla\chi_\eps, ((\nabla\chi_\eps)^\top\nabla\chi_\eps)^{\dot{}},c_\eps\big)\dx\dt\\
			&= \frac{C}{\eps^2}\int_0^T\int_\Omega \zeta(\nabla\chi_\eps,\nabla\dot\chi_\eps,c_\eps)\dx\dt.
		\end{align*}			
The bound now follows from the the energy-dissipation equality \eqref{Eqn:EnDissBal}.

\item[(iv)] Again using the boundedness of $\calE_\eps$, we see that for a.e. $t\in[0,T]$
\[
\norm{\rmD^2\chi_\eps}^p_{L^p(\Omega)} \leq C\int_\Omega \scrH(\rmD^2\chi_\eps)\dx \leq C\eps^2.
\]
Thus, it follows that $\norm{\rmD^2\chi_\eps}_{L^\infty(0,T;L^p(\Omega))}\leq C\eps^\frac{2}{p}$. 
\end{itemize}
\end{proof}

We now prove a priori estimates for the concentration $c_\eps$ and the rescaled concentration variation $\rho_\eps$.

\begin{Lemma}[A priori estimates, part 2]
\label{Lemma:APriori2}
Let $(\chi_\eps, c_\eps)$ be a weak solution of \eqref{Eqn:FS:WeakScaled} in the sense of Definition \ref{Def:weakSolutionFiniteStrain}. Then, there exists a constant $C>0$ (independent of $\eps$) such that:
\begin{itemize}
\item[(i)] $\norm{c_\eps}_{L^\infty(0,T;L^\infty(\Omega))} \leq C$,
\item[(ii)] $\norm{c_\eps\log(\frac{c_\eps}{c_\eq})-c_\eps+c_\eq}_{L^\infty(0,T;L^1(\Omega))} \leq C\eps^2$,
\item[(iii)] $\norm{\rho_\eps}_{L^\infty(0,T;L^2(\Omega))} \leq C$,
\item[(iv)] $\Norm{\nabla c_\eps^{m/2}}_{L^2([0,T]\times\Omega)} + \gamma_1\Norm{\nabla c_\eps^{\frac{m+1+r}{2}}}_{L^2([0,T]\times\Omega)} + \gamma_1\Norm{\nabla c_\eps^{\frac{m}{2}+1+r}}_{L^2([0,T]\times\Omega)} \leq C\eps^{\frac{2}{p}}$,
\item[(v)] $\norm{\calM(\nabla\chi_\eps,c_\eps)\nabla\mu_{*,\eps}}_{L^2([0,T]\times\Omega)} \leq C$,
\item[(vi)] $\norm{\dot \rho_\eps}_{L^{2}(0,T;H^{1}(\Omega)^*)} \leq C$.
\end{itemize}
If, additionally, $\gamma_1>0$ in Assumption \ref{Assu:FreeEn}(ii), then we also have that
\begin{itemize}
\item[(vii)] $\Norm{\frac{1}{r+2}(c_\eps^{r+2}-c_\eq^{r+2})-c_\eq^{r+1}(c_\eps-c_\eq)}_{L^\infty(0,T;L^1(\Omega))} \leq C\eps^2$.
\end{itemize}
\end{Lemma}

\begin{proof}
\begin{itemize}
\item[(i)] This is a direct consequence of Theorem \ref{Th:ConcBounded}.

\item[(ii)] This follows by integrating the lower bound $\pl^2_{cc}\Phi \gtrsim \frac{1}{c}$ from Assumption \ref{Assu:FreeEn}(ii) twice, using that $\Phi(I,c_\eq)=0$, $\pl_c\Phi(I,c_\eq)=0$, the uniform boundedness of $\calE_\eps$, and the \\$L^\infty(0,T;L^2(\Omega;\R^{d\times d}))$-bound for $\nabla u_\eps$.

\item[(iii)] We use the $L^\infty(0,T;L^\infty(\Omega))$-bound for $c_\eps$ from (i) to find a $K=K(\norm{c_\eps}_{L^\infty(0,T;L^\infty(\Omega))})>0$ such that $c_\eps\log(\frac{c_\eps}{c_\eq})-c_\eps+c_\eq \geq K(c_\eps-c_\eq)^2$. The bound now follows from (ii).

\item[(iv)] By the energy-dissipation inequality \eqref{Eqn:EnDissBal}, we have
\[
\int_0^T\int_\Omega \calM(\nabla\chi_\eps,c_\eps)\nabla\mu_\eps\cdot\nabla\mu_\eps \dx\dt \leq C\eps^2.
\]
Using that $\nabla\mu_\eps = \pl^2_{Fc}\Phi(\nabla\chi_\eps,c_\eps)\rmD^2\chi_\eps + \pl^2_{cc}\Phi(\nabla\chi_\eps,c_\eps)\nabla c_\eps$, we have
\begin{align*}
\int_0^T\int_\Omega \calM(\nabla\chi_\eps,c_\eps)\nabla\mu_\eps\cdot\nabla\mu_\eps \dx\dt 
&\geq C\int_0^T\int_\Omega \pl^2_{cc}\Phi(\nabla\chi_\eps,c_\eps)^2 \calM(\nabla\chi_\eps,c_\eps) \abs{\nabla c_\eps}^2\dx\dt \\
&\qquad- C\int_0^T\int_\Omega \abs{\calM(\nabla\chi_\eps,c_\eps)}\Abs{\pl^2_{Fc}\Phi(\nabla\chi_\eps,c_\eps)\rmD^2\chi_\eps}^2 \dx\dt,
\end{align*}
where we have used Young's inequality to absorb the mixed term.

Next, we note that by Assumption \ref{Assu:MobilityTensor} and \ref{Assu:FreeEn}(ii)
\begin{align*}
\int_0^T\int_\Omega \pl^2_{cc}&\Phi(\nabla\chi_\eps,c_\eps)^2 \calM(\nabla\chi_\eps,c_\eps) \abs{\nabla c_\eps}^2\dx\dt \\
&\geq C\big(\Norm{\nabla c_\eps^{\frac{m}{2}}}_{L^2(0,T;L^2(\Omega))}^2 + \gamma\Norm{\nabla c_\eps^{\frac{m+1+r}{2}}}_{L^2(0,T;L^2(\Omega))}^2 + \gamma\Norm{\nabla c_\eps^{\frac{m}{2}+1+r}}_{L^2(0,T;L^2(\Omega))}^2\big).
\end{align*}
Using the $L^\infty(0,T;L^\infty(\Omega))$-bound for $c_\eps$ from (i) and that $m+2\alpha\geq 0$ (Assumption \ref{Assu:FreeEn}(iii), we now obtain
\begin{align*}
\int_0^T\int_\Omega \abs{\calM(\nabla\chi_\eps,c_\eps)}\Abs{\pl^2_{Fc}\Phi(\nabla\chi_\eps,c_\eps)\rmD^2\chi_\eps}^2 \dx\dt &\leq C\norm{c_\eps}^{m+2\alpha}_{L^\infty(0,T;L^{\infty}(\Omega))}\norm{\rmD^2\chi_\eps}^2_{L^\infty(0,T;L^p(\Omega))}\\
&\leq C\eps^{\frac{4}{p}}.
\end{align*}
Combining everything, we have
\begin{align*}
C\eps^\frac{4}{p} &\geq \int_0^T\int_\Omega \calM(\nabla\chi_\eps,c_\eps)\nabla\mu_\eps\cdot\nabla\mu_\eps\dx\dt\\
&\geq C\big(\norm{\nabla c_\eps^{\frac{m}{2}}}_{L^2(0,T;L^2(\Omega))}^2 + \gamma\norm{\nabla c_\eps^{\frac{m+1+r}{2}}}_{L^2(0,T;L^2(\Omega))}^2 + \gamma\norm{\nabla c_\eps^{\frac{m}{2}+1+r}}_{L^2(0,T;L^2(\Omega))}^2\big),
\end{align*}
from which (iv) follows.

\item[(v)] The energy-dissipation inequality \eqref{Eqn:EnDissBal} implies that $\norm{\calM(\nabla\chi_\eps,c_\eps)^\frac{1}{2}\nabla\mu_{*,\eps}}_{L^2([0,T]\times\Omega)}\leq C$. Thus, using the $L^\infty(0,T;L^\infty(\Omega))$-bound for $c_\eps$, it then follows that
\begin{align*}
&\norm{\calM(\nabla\chi_\eps,c_\eps)\nabla\mu_{*,\eps}}_{L^2([0,T]\times\Omega)}
\leq \Norm{c_\eps}^\frac{m}{2}_{L^\infty(0,T;L^\infty(\Omega))}\norm{\calM(\nabla\chi_\eps,c_\eps)^\frac{1}{2}\nabla\mu_{*,\eps}}_{L^2([0,T]\times\Omega)} \leq C.
\end{align*}

\item[(vi)] Testing \eqref{Eqn:FS:Con:Scaled} with $\psi\in L^2(0,T;W^{1,s'}(\Omega))$ we obtain
\begin{align*}
\int_0^T \ip{\dot \rho_\eps,\psi}\dt &\leq C\norm{\calM(\nabla\chi_\eps,c_\eps)\nabla\mu_{*,\eps}}_{L^2([0,T]\times\Omega)}\norm{\psi}_{L^2(0,T;H^{1}(\Omega))} \leq C.
\end{align*}
In fact, we see that $\dot\rho_\eps\in L^2(0,T;H^1(\Omega)^*)$, and thus it suffices to take test functions $\psi\in L^2(0,T;H^1(\Omega))$.

\item[(vii)] This follows in the same way as (ii), but now using the lower bound $\pl^2_{cc}\Phi \gtrsim \gamma_1 c^r$ from Assumption \ref{Assu:FreeEn}(ii). Indeed, $\int_{c_\eq}^c \tilde{c}^r\dd \tilde{c} = \frac{1}{r+1}(c^{r+1}-c_\eq^{r+1})$, and $\int_{c_\eq}^c (\tilde{c}^{r+1}-c_\eq^{r+1})\dd \tilde{c} = \frac{1}{r+2}(c^{r+2}-c_\eq^{r+2}) -c_\eq^{r+1}(c-c_\eq)$.
\end{itemize} 
\end{proof}

\subsection{Limit passage \texorpdfstring{$\eps\to 0$}{eps to 0}}

We are now ready to pass to the limit $\eps\to 0$ in the system \eqref{Eqn:FS:WeakScaled}.

\begin{proof}[Proof of Thm.~\ref{Th:LimitPassage}(i)]

Using the a priori estimates, we can extract converging subsequences (not relabeled) and limits $(u,\rho)$ such that:
\begin{alignat*}{2}
&u_\eps \wstarto u \qquad &&\text{in }L^\infty(0,T;H^1(\Omega)),\\
&\nabla\dot u_\eps \wto \nabla\dot u \qquad &&\text{in }L^2([0,T]\times\Omega),\\
&\rho_\eps \wstarto \rho \qquad &&\text{in }L^\infty(0,T;L^2(\Omega)),\\
&\dot \rho_\eps\wto \dot \rho\qquad &&\text{in }L^2(0,T;H^1(\Omega)^*).%
\end{alignat*}

Note that using the Aubin--Lions lemma %
 we can extract a further subsequence (not relabeled) such that 
\begin{alignat*}{2}
&u_\eps \sto u \qquad &&\text{in }C(0,T;L^2(\Omega)).
\end{alignat*}

\noindent\underline{Mechanical equation.}
We first pass to the limit $\eps\to 0$ in the mechanical equation \eqref{Eqn:FS:Def:Scaled}. 
Using Assumption \ref{Assu:Regul} and Taylor's theorem, we find that for $\C:=\pl^2_{FF}\Phi(I,c_\eq)$ and $\K:=\pl^2_{Fc}\Phi(I,c_\eq)$
\begin{align}
\label{Eqn:Taylor:Interm1}
\abs{\pl_F\Phi(\nabla\chi_\eps,c_\eps) - \eps\C\nabla u_\eps - \eps\K\rho_\eps} \leq C\eps^2(\abs{\nabla u_\eps}^2 + \abs{\rho_\eps}^2).
\end{align}
In particular, using the rigidity estimate in Lemma \ref{Lemma:RigidityEst}(ii) and the $L^\infty(0,T;L^\infty(\Omega))$-bound for $c_\eps$, it follows that 
\begin{align*}
\abs*{\frac{1}{\eps}\pl_F\Phi(\nabla\chi_\eps,c_\eps) - \C\nabla u_\eps - \K\rho_\eps} \leq C(\abs{\nabla u_\eps} + \abs{\rho_\eps}),
\end{align*}
and thus $\frac{1}{\eps}\pl_F\Phi(\nabla\chi_\eps,c_\eps) - \C\nabla u_\eps - \K\rho_\eps$ is bounded in $L^\infty(0,T;L^2(\Omega))$ and has a weak*-converging subsequence. However, \eqref{Eqn:Taylor:Interm1} also implies that
\begin{align*}
\norm*{\frac{1}{\eps}\pl_F\Phi(\nabla\chi_\eps,c_\eps) - \C\nabla u_\eps - \K\rho_\eps}_{L^\infty(0,T;L^1(\Omega))} &= \sup_{t\in[0,T]}\int_\Omega  \abs*{\frac{1}{\eps}\pl_F\Phi(\nabla\chi_\eps,c_\eps) - \C\nabla u_\eps - \K\rho_\eps}\dx \\
&\leq C\eps(\norm{\nabla u_\eps}^2_{L^\infty(0,T;L^2(\Omega))} + \norm{\rho_\eps}^2_{L^\infty(0,T;L^2(\Omega))})\\
&\leq C\eps\to 0,
\end{align*}
and thus $\frac{1}{\eps}\pl_F\Phi(\nabla\chi_\eps,c_\eps) - \C\nabla u_\eps - \K\rho_\eps$ must weak*-converge to 0 in $L^\infty(0,T;L^2(\Omega))$.
Using the weak*-convergence of $u_\eps$ in $L^\infty(0,T;H^1(\Omega))$ and the weak*-convergence of $\rho_\eps$ in $L^\infty(0,T;L^2(\Omega))$, we then conclude that 
\begin{align*}
\int_0^T\int_\Omega \frac{1}{\eps}\pl_F\Phi(\nabla\chi_\eps):\nabla\phi\dx\dt \to \int_0^T\int_\Omega \big(\C \nabla u + \K\rho\big) : \nabla\phi \dx\dt. 
\end{align*}

The limit passage in the viscous stress follows in a similar way, now using the weak convergence of $\nabla\dot u_\eps$ in $L^2([0,T]\times\Omega)$. Setting $\D := \pl^2_{\dot{F}\dot{F}}\zeta(I,0,c_\eq)$ (which is related to $\wt\D$ in Assumption \ref{Assu:ViscousStress} by $\D=4\wt\D$), it follows that 
\[
\int_0^T\int_\Omega \frac{1}{\eps}\pl_{\dot F}\zeta(\nabla\chi_\eps,\nabla\dot\chi_\eps,c_\eps):\nabla\phi\dx\dt \to \int_0^T\int_\Omega \D\nabla\dot u:\nabla\phi \dx\dt.
\]

To show that the hyperstress vanishes as $\eps\to 0$, we note that H\"older's inequality and the a priori estimate for $\rmD^2u_\eps$ imply that
\begin{align*}
\int_0^T\int_\Omega \frac{1}{\eps}\pl_G\scrH(\rmD^2\chi_\eps)\vdots \rmD^2\phi\dx\dt &\leq \frac{1}{\eps}\int_0^T \norm{\pl_G\scrH(\rmD^2\chi_\eps)}_{L^\frac{p}{p-1}(\Omega)}\norm{\rmD^2\phi}_{L^p(\Omega)}\dt\\
&\leq \frac{C}{\eps}\int_0^T\norm{\rmD^2\chi_\eps}_{L^p(\Omega)}^{p-1}\norm{\phi}_{W^{2,p}(\Omega)}\dt\\
&\leq C\eps^{-1+\frac{2(p-1)}{p}}\\
&\leq C\eps^{1-\frac{2}{p}}\to 0.
\end{align*}

This concludes the limit passage for the mechanical equation. 
\medskip 

\noindent\underline{Diffusion equation.}
We now pass to the limit in the diffusion equation \eqref{Eqn:FS:Con:Scaled}, i.e., in 
\[
\int_0^T \ip{\dot \rho_\eps,\psi}\dt + \int_0^T\int_\Omega\calM\big(I+\eps\nabla u_\eps,c_\eq+\eps\rho_\eps\big)\nabla\mu_{*,\eps} \cdot \nabla \psi \dx\dt = 0.
\] 
The limit passage in the first integral follows directly from the weak convergence of $\rho_\eps$ in $L^2(0,T;H^1(\Omega)^*)$. 
To pass to the limit in the second integral, we use Lemma \ref{Lemma:APriori2}(v) to find $\xi\in L^2([0,T]\times\Omega)$ such that
\[
\calM(\nabla\chi_\eps,c_\eps)\nabla\mu_{*,\eps} \wto \xi \quad \text{in }L^2([0,T]\times\Omega).
\]
We now identify the limit $\xi$. This is done in three steps.

\underline{Case I: $\gamma_1=\gamma_2=0$:}

\textit{Step 1.} Using the bounds for $\Phi$ and the fact that $\pl_{c}\Phi(I,c_\eq)=0$ (see Assumption \ref{Assu:FreeEnEq}), we see that
\begin{align*}
\mu_\eps = \pl_{c}\Phi(\nabla\chi_\eps,c_\eps) &= \pl_{c}\Phi(\nabla\chi_\eps,c_\eps) - \pl_{c}\Phi(I,c_\eps) + \pl_{c}\Phi(I,c_\eps) - \pl_{c}\Phi(I,c_\eq)\\
&\leq C\Big(c_\eps^\alpha\abs{\nabla\chi_\eps - I} + \log\frac{c_\eps}{c_\eq}\Big) .
\end{align*}
Thus, using the $L^\infty(0,T;L^\infty(\Omega;\R^{d\times d})$-bound for $\nabla\chi_\eps$ and Assumption \ref{Assu:MobilityTensor}, we see that
\begin{align*}
\norm{\calM(\nabla\chi_\eps,c_\eps)\mu_\eps}^2_{L^2([0,T]\times\Omega)} &\leq C\int_0^T\int_\Omega c_\eps^{2m}\Abs{\log\frac{c_\eps}{c_\eq}}^2\dx\dt + C\eps^2\int_0^T\int_\Omega c_\eps^{2(m+\alpha)}\abs{\nabla u_\eps}^2\dx\dt\\
&\leq C\norm{c_\eps}^{m}_{L^\infty(0,T;L^\infty(\Omega))}\int_0^T\int_\Omega c_\eps^{m}\Abs{\log\frac{c_\eps}{c_\eq}}^2\dx\dt \\
&\qquad+ C\eps^2\norm{c_\eps}^{2(m+\alpha)}_{L^\infty(0,T;L^\infty(\Omega))}\norm{u_\eps}^2_{L^2(0,T;H^1(\Omega))}\\
&\leq C\norm{c_\eps-c_\eq}^2_{L^2([0,T]\times\Omega)} + C\eps^2 \leq C\eps^2,
\end{align*}
where we have used Lemma \ref{Lemma:Appendix1}, i.e., that for $m\leq 2-\eta$ there exists some $C>0$ such that $c_\eps^m\abs{\log\frac{c_\eps}{c_\eq}}^2\leq C\abs{c_\eps-c_\eq}^2$.

In particular, we see that $\frac{1}{\eps}\calM(\nabla\chi_\eps,c_\eps)\mu_\eps$ has a weakly converging subsequence in $L^2([0,T]\times\Omega)$. We now define $\mu_*$ as the weak limit 
\[
\M(I,c_\eq)^{-1}\calM(\nabla\chi_\eps,c_\eps)\frac{\mu_\eps}{\eps}\wto \mu_* \quad\text{ in }L^2([0,T]\times\Omega).
\]
Since $\calM$ is $C^1$, we can use Taylor's theorem and $\mu_\eps = \pl_c\Phi(\nabla\chi_\eps,c_\eps)$ to find
\[
\Abs{\M(I,c_\eq)^{-1}\calM(\nabla\chi_\eps,c_\eps)\frac{\mu_\eps}{\eps} - \K\nabla u_\eps - \LL\rho_\eps} \leq C\eps(\abs{\nabla u_\eps}^2+\abs{\rho_\eps}^2).
\]
Using the $L^\infty(0,T;L^2(\Omega;\R^{d\times d}))$-bound for $\nabla u_\eps$ and the $L^\infty(0,T;L^2(\Omega))$-bound for $\rho_\eps$, it then follows that $\mu_* = \K\nabla u + \LL\rho$ for a.e. $(t,x)\in [0,T]\times\Omega$.

\textit{Step 2.} Let $1<s<\min\{\frac{2p}{p+2},\frac{p}{p-1}\}<2$. We will show that
\[
\calM(\nabla\chi_\eps,c_\eps)\nabla\mu_\eps \wto \M(I,c_\eq)\nabla\mu_* \quad \text{in }L^s([0,T]\times\Omega).
\]
Then, uniqueness of weak limits implies that $\xi = \M(I,c_\eq)\nabla\mu_*$.

We denote by $D$ the derivative for which $D\calM(\nabla\chi,c) = D_F\calM(\nabla\chi,c)\rmD^2\chi + D_c\calM(\nabla\chi,c)\nabla c$. Using the bounds in \ref{Assu:MobReg}, we then have
\begin{align*}
\norm{\rmD(\calM(\nabla\chi_\eps,c_\eps))\mu_\eps}_{L^s([0,T]\times\Omega)} &\leq C\Norm{c_\eps^m\abs{\rmD^2\chi_\eps}\mu_\eps}_{L^s([0,T]\times\Omega)} + C\Norm{c_\eps^{m-1}\abs{\nabla c_\eps}\mu_\eps}_{L^s([0,T]\times\Omega)}\\
&=: I_1 + I_2.
\end{align*}
We estimate the first integral $I_1$. Since $s<\frac{2p}{p+2}$, we use the embedding $L^\frac{2p}{p+2}(\Omega)\hookrightarrow L^s(\Omega)$, H\"older's inequality, and the estimate $\norm{\rmD^2\chi_\eps}_{L^\infty(0,T;L^p(\Omega;\R^{d\times d\times d}))}\leq C\eps^{\frac{2}{p}}$ from Lemma \ref{Lemma:APriori1}(iv) to obtain
\begin{align*}
I_1 &\leq C\Norm{c_\eps^m\abs{\rmD^2\chi_\eps}\mu_\eps}_{L^\frac{2p}{p+2}([0,T]\times\Omega)}\\
&\leq C\norm{\rmD^2\chi_\eps}_{L^\infty(0,T;L^p(\Omega;\R^{d\times d\times d}))}\Norm{c_\eps}^{\frac{m}{2}}_{L^\infty(0,T;L^\infty(\Omega))}\Norm{c_\eps^\frac{m}{2}\mu_\eps}_{L^2([0,T]\times\Omega)}\\
&\leq C\eps^{\frac{2}{p}}\Norm{c_\eps^\frac{m}{2}\mu_\eps}_{L^2([0,T]\times\Omega)}.
\end{align*}
We again note that $\mu_\eps \leq  C\big(c_\eps^\alpha\abs{\nabla\chi_\eps - I} + \log\frac{c_\eps}{c_\eq}\big)$, and proceed as in Step 1 to obtain \\$\norm{c_\eps^\frac{m}{2}\mu_\eps}_{L^2([0,T]\times\Omega)} \leq C\eps$. In particular, we see that $I_1\leq C\eps^{1+\frac{2}{p}}$.

To estimate the second integral $I_2$, we note that $1<s<\frac{p}{p-1}<2$, and thus
\begin{align*}
I_2 &\leq C\Norm{c_\eps^{m-1}\abs{\nabla c_\eps}\mu_\eps}_{L^\frac{p}{p-1}([0,T]\times\Omega)}\leq C\Norm{\nabla c_\eps^{\frac{m}{2}}}_{L^2([0,T]\times\Omega)}\Norm{c_\eps^\frac{m}{2}\mu_\eps}_{L^q([0,T]\times\Omega)},
\end{align*}
where $2<q<\frac{2p}{p-2}$. We now note that
\begin{align*}
\Norm{c_\eps^\frac{m}{2}\mu_\eps}_{L^q([0,T]\times\Omega)} &= \Norm{c_\eps^\frac{m}{2}\mu_\eps}^{1-\frac{2}{q}}_{L^\infty(0,T;L^\infty(\Omega))}\Norm{c_\eps^\frac{m}{2}\mu_\eps}^{\frac{2}{q}}_{L^2([0,T]\times\Omega)} \leq C\eps^\frac{2}{q},
\end{align*}
where the $L^\infty(0,T;L^\infty(\Omega))$-bound for $c_\eps^\frac{m}{2}\mu_\eps$ follows from the $L^\infty(0,T;L^\infty(\Omega))$-bound for $c_\eps$ and the $L^\infty(0,T;L^\infty(\Omega;\R^{d\times d})$-bound for $\nabla\chi_\eps-I$ (see Lemma \ref{Lemma:RigidityEst}(ii)). In particular, using the estimate in Lemma \ref{Lemma:APriori2}(iv), we obtain
\[
I_2 \leq C\eps^{\frac{2}{p}+\frac{2}{q}} = C\eps^{d},
\]
where $d=\frac{2}{p}+\frac{2}{q}>\frac{2}{p}+2\frac{p-2}{2p}=1$.

Combining everything, we thus see that for some $\wt d>1$.
\[
\norm{\rmD(\calM(\nabla\chi_\eps,c_\eps))\mu_\eps}_{L^s([0,T]\times\Omega)} \leq C\eps^{\wt d},
\]
which implies that
\[
\frac{1}{\eps}\rmD(\calM(\nabla\chi_\eps,c_\eps))\mu_\eps \sto 0\quad \text{in } L^s([0,T]\times\Omega).
\]

\textit{Step 3.}
Denote by $\rmD(\calM\mu)\in \R^{d\times d\times d}$ the derivative defined by $(\rmD(\calM\mu))_{ijk} = \pl_k(\calM_{ij}\mu)$. The product rule then implies that
\[
\rmD(\calM(\nabla\chi_\eps,c_\eps)\mu_\eps) = \rmD(\calM(\nabla\chi_\eps,c_\eps))\mu_\eps + \calM(\nabla\chi_\eps,c_\eps)\nabla\mu_\eps.
\]
Using the estimate from Step 2 for the first term, and Lemma \ref{Lemma:APriori2}(v) for the second term, it then follows that
\[
\norm{\calM(\nabla\chi_\eps,c_\eps)\mu_\eps}_{L^s(0,T;W^{1,s}(\Omega;\R^{d\times d}))} \leq C\eps,
\]
and thus $\frac{1}{\eps}\calM(\nabla\chi_\eps,c_\eps)\mu_\eps$ has a weakly convergent subsequence in $L^s(0,T;W^{1,s}(\Omega;\R^{d\times d}))$. From Step 1 we know that the $L^2([0,T]\times\Omega)$-limit is $\M(I,c_\eq)\mu_*$, and thus it follows that
\[
\frac{1}{\eps}\calM(\nabla\chi_\eps,c_\eps)\mu_\eps \wto \M(I,c_\eq)\mu_* \quad \text{in }L^s(0,T;W^{1,s}(\Omega)).
\]
Finally, we conclude that
\begin{align*}
\frac{1}{\eps}\calM(\nabla\chi_\eps,c_\eps)\nabla\mu_\eps &= \frac{1}{\eps}\rmD(\calM(\nabla\chi_\eps,c_\eps)\mu_\eps) - \frac{1}{\eps}\rmD(\calM(\nabla\chi_\eps,c_\eps))\mu_\eps \\
&\wto \M(I,c_\eq)\nabla\mu_* \quad \text{in }L^s([0,T]\times\Omega),
\end{align*}
and thus $\xi = \M(I,c_\eq)\nabla\mu_*$.
This concludes the limit passage in the diffusion equation \eqref{Eqn:FS:Con:Scaled} for the case that $\gamma_1=\gamma_2=0$.

\underline{Case II: $\gamma_2\geq\gamma_1>0$:} We now highlight the changes for the case that $\gamma_2\geq\gamma_1>0$. Now, we have that
\begin{align*}
\mu_\eps &\leq C\Big(c_\eps^\alpha\abs{\nabla\chi_\eps - I} + \log\frac{c_\eps}{c_\eq} + (c_\eps^{r+1}-c_\eq^{r+1})\Big).
\end{align*}
The third term, however, can be estimated using Lemma \ref{Lemma:Appendix2}, i.e., there exists a constant $C>0$ such that
\[
\abs{c_\eps^{r+1}-c_\eq^{r+1}} \leq C\Big( \Big|\frac{c_\eps^{r+2}-c_\eq^{r+2}}{r+2} - c_\eq^{r+1}(c_\eps-c_\eq)\Big| + \abs{c_\eps - c_\eq} \Big).
\]
The proof now follows in the same way as in Case I, additionally using the bound in Lemma \ref{Lemma:APriori2}(vii).
\end{proof}

\subsection{Properties of \texorpdfstring{$\C$, $\K$ and $\D$}{tensors}}

We now show some properties of the tensors $\C$, $\K$ and $\D$. In particular, we prove that these tensors only act on the symmetric part of matrices, and that $\C$ is positive definite. For a matrix $U\in\R^{d\times d}$, we denote by $U^{\sym}:=\frac{1}{2}(U+U^\top)$ the symmetric part of $U$, and by $U^{\an}:=\frac{1}{2}(U-U^\top)$ the antisymmetric part. 

\begin{Lemma}
\label{Lemma:TensorsActOnSymm}
Let $U\in\GL^+(d)$. Then, the tensors $\C$, $\K$ and $\D$ only act on the symmetric part of $U$, i.e., $\C U = \C U^\sym$, etc. %
\end{Lemma}
\begin{proof}
We modify the derivation leading up to \cite[Eqn.~(2.3)]{MieSte2013LPitELoFP}. For example, to prove that $\C U^{\text{anti}}=0$, note that for any $\lambda\in \R$ we have $\exp(\lambda U^\an)\in \SO(d)$. %
So, letting $\xi(\lambda) := \pl_F\Phi(\exp(\lambda U^\an),c_\eq)$, the frame indifference of $\Phi$ and $\pl_F\Phi(I,c_\eq)=0$ imply that $\xi \equiv 0$. In particular, differentiating with respect to $\lambda$ and setting $\lambda=0$ gives $\C U^\an = 0$. The proofs for the other tensors follow in a similar way.
\end{proof}

\begin{Lemma}[Positive definiteness of $\C$]\label{Lemma:PropC}
Let $U\in\GL^+(d)$. Then, there exists a $C>0$ such that $\C U^\sym:U^\sym\geq C\abs{U^\sym}^2$.
\end{Lemma}
\begin{proof}
We follow \cite[Eqn.~(2.4)-(2.5)]{MieSte2013LPitELoFP}. Note that by Lemma \ref{Lemma:TensorsActOnSymm}, $\C$ only acts on $U^\sym$, i.e., $\C U=\C U^\sym$.

By linearizing the distance function $\dist$ around the identity (see e.g. \cite[Eqn.~3.20]{FrRiMu2002GRDN}), we obtain 
\[
\frac{1}{\eps}\dist(I+\eps U,\SO(d)) = U^\sym + \eps\calO(\abs{U}^2).
\]
Next, using that $\Phi$ is $C^3$ in a neighbourhood of $I\times c_\eq$ it follows that 
\[
\frac{1}{\eps^2}\Phi(I+\eps U,c_\eq) = \frac{1}{2}\C U:U + \eps\calO(\abs{U}^3).
\]
Using these two equalities and Assumption \ref{Assu:PhiNondegen}, it follows that
\begin{align*}
C\abs{U^\sym}^2 &\leq \lim_{\eps\to 0} \frac{C}{\eps^2} \dist^2(I+\eps U,\SO(d)) \leq \lim_{\eps\to 0}\frac{1}{\eps^2} \Phi(I+\eps U,c_\eq) = \frac{1}{2}\C U^\sym:U^\sym,
\end{align*}
completing the proof.
\end{proof}

\subsection{Uniqueness of small-strain solutions}

It remains to show that the obtained weak solution $(u,\rho)$ of the linearized problem \eqref{Eqn:SS:WeakSoln} is the unique weak solution. To do this, we start by proving that weak solutions of the linear system \eqref{Eqn:SS} satisfy an energy-dissipation balance.

\begin{Lemma}[Energy-dissipation balance]
Let $(u,\rho)$ be a weak solution of the linearized problem \eqref{Eqn:SS:WeakSoln} in the sense of Definition \ref{Def:WeakSolnLin}. Then, the following energy balance is satisfied:
\begin{multline}
\label{Eqn:EnBal:Lin}
\calE_0(t,u,\rho) + 2\int_0^t \calR_0(\dot u)\ds + \int_0^t\int_\Omega\M(I,c_\eq)\nabla \mu_*\cdot \nabla\mu_*\dx\ds 
= \calE_0(0,u_0,\rho_0) - \int_0^t \ip{\dot\ell_*(s),u}\ds,
\end{multline}
where $\calE_0$ and $\calR_0$ are the quadratic forms introduced in \eqref{Eqn:E_0} and \eqref{Eqn:R_0}, i.e.,
\begin{align*}
\calE_0(t,u,\rho) &= \int_\Omega \frac{1}{2}\C e(u):e(u) + \rho\K:e(u) + \frac{1}{2}\LL \rho^2 \dx - \ip{\ell_*(t),u},\\
\calR_0(\dot u) &= \int_\Omega \frac{1}{2}\D e(\dot u):e(\dot u)\dx. 
\end{align*}
\end{Lemma}
\begin{proof}
A direct calculation shows that
\[
\frac{\dd}{\dt}\calE_0(t,u,\rho) = \int_\Omega \C e(u):e(\dot u) + \rho\K: e(\dot u)\dx + \ip{\dot\rho,\K e(u) + \LL\rho} - \ip{\dot\ell_*(t),u} - \ip{\ell_*(t),\dot u}. 
\]
Integrating from $s=0$ to $s=t$ and using that $\dot u$ and $\mu_* = \K{:} e(u) + \LL\rho$ are valid test functions for \eqref{Eqn:SS:WeakDef} and \eqref{Eqn:SS:WeakCon}, respectively, the energy-dissipation balance \eqref{Eqn:EnBal:Lin} follows.
\end{proof}

Using this energy-dissipation balance, we can now show that there exists at most one weak solution.

\begin{proof}[Proof of Thm.~\ref{Th:LimitPassage}(ii)]
Suppose that there exist two weak solutions $(u_1,\rho_1)$ and $(u_2,\rho_2)$. 
Define $u_3:=u_1-u_2$, $\rho_3:=\rho_1-\rho_2$, then $(u_3,\rho_3)$ is a solution of
\[
\int_0^T\int_\Omega (\C e(u) + \K\rho + \D e(\dot u)):\nabla\phi \dx\dt = 0,
\]
\[
\int_0^T \ip{\dot\rho,\psi}\dt + \int_0^T\int_\Omega \M(I,c_\eq)\nabla\mu_*(u,\rho)\cdot\nabla\psi\dx\dt %
= 0.
\]
Thus, for 
\[
\mathscr{E}(u,\rho) %
= \int_\Omega \frac{1}{2}\C e(u):e(u) + \K e(u)\rho + \frac{1}{2}\LL\rho^2 \dx \geq 0
\]
it follows that for almost every $t\in [0,T]$:
\begin{align*}
\mathscr{E}(u_3(t),\rho_3(t)) &- \mathscr{E}(u_3(0),\rho_3(0)) =
\int_0^t\frac{\dd}{\ds}\mathscr{E}(u_3(s),\rho_3(s))\ds \\
&= -\int_0^t\int_\Omega \D e(\dot u_3):e(\dot u_3)\dx\ds - \int_0^t\int_\Omega \M(I,c_\eq)\abs{\nabla\mu_*(u_3,\rho_3)}^2\dx\ds \\
&\leq 0.
\end{align*}
In particular, as $\mathscr{E}(u_3(0),\rho_3(0)) = \mathscr{E}(0,0) = 0$ and $\mathscr{E}$ is nonnegative, it follows that $\mathscr{E} \equiv 0$, i.e., $u_3 \equiv 0$, $\rho_3\equiv 0$.

\end{proof}

As a consequence of the energy inequality, we also note that solutions of the linearized system \eqref{Eqn:SS} converge to an equilibrium state whenever $t\to \infty$ (if solutions exist for all times $t\geq 0$), and the boundary data are chosen in a suitable way.

\begin{Cor}
Suppose that the external forces $f_*$ and $g_*$ are time-independent. Let $(u,\rho)$ be a weak solution of \eqref{Eqn:SS:WeakSoln} in the sense of Definition \ref{Def:WeakSolnLin} that exists for all times $t\geq 0$. Then, $(u,\rho)$ converges to a solution $(v,\xi)$ of the (static) system:
\begin{subequations}
\begin{alignat*}{2} 
- \DIV\!\big(\C e(v) + \K\xi \big) &= f_*\quad &&\text{in}\ \Omega,\\
-\DIV\!\big(\M(I,c_\eq)\nabla \nu_*\big) &= 0\quad &&\text{in}\ \Omega,
\end{alignat*}
\end{subequations}
where $\nu_* = \K:e(v)+\LL\xi$.
\end{Cor}
\begin{proof}
Define $w(t)=u(t)-v$, $\zeta(t)=\rho(t)-\xi$, $\lambda_*(t)=\mu_*(t)-\nu_*$ then $(w,\zeta)$ satisfies the system
\[
\int_0^T\int_\Omega (\C e(w) + \K\zeta + \D e(\dot w)):\nabla\phi \dx\dt = 0
\]
\[
\int_0^T \ip{\dot\zeta,\psi}\dt + \int_0^T\int_\Omega \M(I,c_\eq)\nabla\lambda_*(w,\zeta)\cdot\nabla\psi\dx\dt = 0. 
\]
In particular, we see that for $t_2>t_1$,
\begin{align*}
\mathscr{E}(w(t_2),\zeta(t_2)) &- \mathscr{E}(w(t_1),\zeta(t_1)) =\int_{t_1}^{t_2}\frac{\dd}{\ds}\mathscr{E}(w(s),\zeta(s))\ds \\
&= -\int_{t_1}^{t_2}\int_\Omega \D e(\dot w):e(\dot w)\dx\ds - \int_{t_1}^{t_2}\int_\Omega \M(I,c_\eq)\abs{\nabla\lambda_*(w,\zeta)}^2\dx\ds \\
&\leq 0.
\end{align*}
Thus, $\mathscr{E}$ is decreasing with time, and since $\mathscr{E}$ is bounded from below by $0$, we see that $\mathscr{E}(t)\to \inf_{t\geq 0}\mathscr{E}(t)$ as $t\to\infty$. However, since $\mathscr{E}(0,0) = 0$, this infimum is attained, and we have $\mathscr{E}(t)\to 0$ as $t\to\infty$. Thus, we have $w(t)\to 0$, $\zeta(t) \to 0$ as $t\to\infty$. 
\end{proof}

\section{Appendix}

\begin{Lemma}
\label{Lemma:Appendix1}
Let $0<m\leq 2-\eta$ for some $\eta>0$, and $c_\eq>0$. Then, there exists a $C>0$ such that for all $x>0$,
\[
x^m\left|\log\frac{x}{c_\eq}\right|^2 \leq C\abs{x-c_\eq}^2\quad \text{for any }x>0.
\]
\end{Lemma}
\begin{proof}
Note that $\lim_{x\to 0}x^m\abs{\log\frac{x}{c_\eq}}^2 = 0$, and thus the desired inequality holds for $x$ in some small interval $(0,\eps)$. To show that the inequality holds for all $x>0$, we restrict to $x\geq \eps$ and define $f(x):=x^\frac{m}{2}\log\frac{x}{c_\eq}$. Since $f'(x)=x^{\frac{m}{2}-1}(\frac{m}{2}\log\frac{x}{c_\eq} + 1)$, it follows that $\wt{C} := \sup_{x\geq\eps}\abs{f'(x)}<\infty$. In particular, since $f$ is continuously differentiable on $[\eps,\infty)$, it follows that $f$ is Lipschitz continuous in $x=c_\eq$ with Lipschitz constant $\wt{C}$. As $f(c_\eq)=0$, we thus obtain
\[
x^\frac{m}{2}\left|\log\frac{x}{c_\eq}\right| = \abs{f(x)-f(c_\eq)} \leq \wt{C}\abs{x-c_\eq}.
\]
Squaring this inequality gives the desired bound.
\end{proof}

\begin{Lemma}
\label{Lemma:Appendix2}
Let $c_\eq>0$ and $r>-1$. Then, there exists a $C>0$ such that for all $x>0$,
\[
\abs{x^{r+1}-c_\eq^{r+1}} \leq C\Big( \Big|\frac{x^{r+2}-c_\eq^{r+2}}{r+2} - c_\eq^{r+1}(x-c_\eq)\Big| + \abs{x - c_\eq} \Big).
\]
\end{Lemma}
\begin{proof}
Using Young's inequality with exponents $\frac{r+1}{r}$ and $r+1$, it follows that
\begin{align*}
x^{r+1}-c_\eq^{r+1} = (r+1)\int_{c_\eq}^x y^{r}\dd y 
&\leq \int_{c_\eq}^x r y^{r+1} + 1\dd y 
= r\frac{x^{r+2}-c_\eq^{r+2}}{r+2} + (x-c_\eq)\\
&= r \Big(\frac{x^{r+2}-c_\eq^{r+2}}{r+2} - c_\eq^{r+1}(x-c_\eq) \Big) + r(c_\eq^{r+1}+1)(x-c_\eq).
\end{align*}
The desired inequality now follows by taking the absolute value and using the triangle inequality.
\end{proof}

\section*{Acknowledgements}
W.v.O. would like to thank Matthias Liero and Annegret Glitzky for helpful discussions.
W.v.O. was partially supported by DFG via the Priority Program SPP 2256 \emph{Variational Methods for Predicting Complex
 Phenomena in Engineering Structures and Materials} (project no. 441470105, subproject Mi 459/9-1 \emph{Analysis for thermomechanical models with internal variables}).

\newcommand{\etalchar}[1]{$^{#1}$}

\end{document}